\documentclass{amsart}
\usepackage{amssymb,amsmath}
\usepackage[american,french]{babel}
\usepackage{amsopn}
\usepackage{fancyhdr}
\title{Ergodicity for  the weakly damped stochastic non-linear Schr\"odinger equations}
\author{Arnaud DEBUSSCHE \quad and  \quad
Cyril ODASSO
\\
 \\  Ecole Normale Sup\'erieure de Cachan, antenne de Bretagne,\\ Avenue Robert Schuman,
 Campus de Ker Lann, 35170 Bruz (FRANCE). \\ and
\\ IRMAR,  UMR 6625 du CNRS, Campus de Beaulieu,  35042 Rennes cedex (FRANCE)}

\newtheorem{Theorem}{Theorem}[section]

\newtheorem{Proposition}[Theorem]{Proposition}
\newtheorem{Lemma}[Theorem]{Lemma}
\newtheorem{Corollary}[Theorem]{Corollary}
\newtheorem{Remark}[Theorem]{Remark}

\makeatletter
\newif\ifmsbmloaded@

\@addtoreset{equation}{section}

\makeatother

\def\R{\mathbb R}
\def\N{\mathbb N}

\def\C{\mathbb C}
\def\E{\mathbb E}
\def\P{\mathbb P}
\def\Pcal{\mathcal{P}}
\def\Hcal{\mathcal{H}}

\def\Dr{\mathcal D}
\def\Br{\mathcal B}

\def\F{\mathcal F}

\def\ds{\displaystyle}

\newcommand{\BLANC}[1]{   }

\newcommand{\abs}[1]{\left\vert#1\right\vert}
\newcommand{\norm}[1]{\left\Vert#1\right\Vert}

\newcommand{\sig}{\sigma}
\renewcommand{\i}{\textrm{i}}
\def \Espace{\renewcommand{\arraystretch}{1.7} }

\newcommand{\abssig}[1]{\abs{#1}_{4}^{4}}

\newcommand{\B}[1]{ \ensuremath{\abs{#1}^{2}  #1}}
\newcommand{\Bd}[1]{ \ensuremath{\abs{#1}^{2}  (#1)}}

\newcommand{\LS}[1]{ \ensuremath{ \frac{d#1}{dt} +\alpha #1+ \i A#1}}
\newcommand{\LSbis}[1]{ \ensuremath{ d#1 +\alpha #1 \;dt+ \i A#1 \;dt}}
\newcommand{\LSabs}[1]{ \ensuremath{ \frac{d\abs{#1}^2}{dt} %
 + 2\alpha\abs{#1}^2}}

\newcommand{\NLS}[1]{ \ensuremath{ \LSbis{#1}-\i  \B{#1}\; dt }}

\newcommand{\NLSDelta}[1]{\ensuremath{d#1 +\alpha #1 \;dt-\i\Delta#1 \;dt%
-\i  \B{#1}\;dt }}

\begin{document}
\setcounter{page}{127}

\selectlanguage{american} \maketitle \pagestyle{fancy}





\noindent\textbf{Abstract}:
 We study a damped stochastic non-linear Schr\"odinger (NLS) equation driven by an additive noise.
 It is white in time and smooth in space.
 Using a coupling method, we establish convergence of the Markov transition semi-group toward a unique invariant
 probability measure. This kind of method was originally developed to prove exponential mixing for strongly dissipative equations such as the Navier-Stokes equations. We consider here a
 weakly dissipative equation, the damped nonlinear Schr\"odinger equation in the
 one-dimensional cubic case.
We prove that the mixing property holds and that the rate of convergence to equilibrium is at least
polynomial of  any power.

\


\section*{Introduction}

The non-linear Schr\"odinger (NLS)
 equation models the propagation of non-linear dispersive waves in various areas of physics such as hydrodynamics
\cite{Newel1}, \cite{Newel2}, optics, plasma physics, chemical reaction \cite{Huber}...

When studying the propagation in random media, a noise can be
introduced. For instance in \cite{FKLT}, \cite{FKLT2}, the cubic
nonlinear Schr\"odinger equation with additive white noise and
damping is introduced. There, the propagation of waves over very
long distance is studied. Damping effect cannot be neglected in
this case and has to be counterbalanced by amplifiers. The white
noise is a model for the description of the randomness in these
amplifiers. Such model is valid if the distance between amplifiers
is small compared to propagation distance.

Our aim in this work is to study ergodicity for this type of equation. We consider the one-dimensional
case with cubic focusing nonlinearity.
It has the form
\begin{equation}\label{EqIntro}
\left \{
\begin{array}{rcll}
\NLSDelta{u} & = & b dW ,& \\
u(t,x) & = & 0,           &\mbox{   for   } x \in \{0,1\}, \; t>0, \\
u(0,x) & = & u_0(x), &\mbox{   for   } x \in [0,1],
\end{array}
\right .
\end{equation}
where $\alpha>0$.
The unknown $u$ is a complex valued process depending on $x \in [0,1]$ and
$t\geq 0$. Dirichlet boundary conditions are considered but we could also use Neumann or
periodic boundary condition.

Existence and uniqueness of solutions for (\ref{EqIntro}) is not very difficult to prove using
straightforward generalization of deterministic arguments.
Note that the damping term is necessary to have an invariant measure. Indeed, if
$\alpha= 0$ and $b\not = 0$ then  the $L^2\left( 0,1 \right)$ norm grows linearly in time.

The Complex Ginzburg-Landau (CGL) is also a form of dissipative NLS equation.
The exponential mixing of the stochastic CGL equation has been establish
in \cite{H} in a particular case and in
\cite{ODASSO} in the general case. The method was inspired by the so-called coupling method.
This method has been introduced in  \cite{BKL},   \cite{H}, \cite{KS},
 \cite{KS2}, \cite{KS3},  \cite{Matt} and \cite{S}. In all these articles, a
strongly dissipative stochastic partial differential equations driven by a noise which may be degenerate
is considered.  Due to the possible degeneracy of the noise Doob Theorem  cannot to be applied
(see \cite{DPZ2} for the theory of ergodicity when Doob Theorem can be applied). Indeed, it requires
the strong Feller property, which can be proved only when the noise lives in a space of
spatially irregular functions, which is clearly not true for a degenerate noise.
The main idea is to
compensate the degeneracy of the noise by dissipativity arguments,
 the so-called  Foias-Prodi estimates. Roughly speaking, the process can be decomposed
 into the sum of a strongly dissipative process and another one driven by a non-degenerate
 noise. The strongly dissipative part is treated as in \cite{DPZ1} section 11.5, while the non-degenerate part is treated thanks to probabilistic arguments. The difficulty is of course in
 the fact that the two parts of the process do not evolve independently so that the two methods
 have to be used simultaneously. The probabilistic part can be treated either by a generalization
 of Doob Theorem (see  \cite{EMS}, \cite{HM04}, \cite{KS00}) or by coupling argument
 (see \cite{KS}, \cite{KS2}, \cite{KS3},  \cite{Matt}). Each method has its advantages. The
 first one allows treating very degenerate noises while the coupling method proves also
 exponential convergence to equilibrium.

In the case of the NLS equation, it seems hopeless to use Doob Theorem. Indeed, due to the lack
of smoothing effect of the deterministic part of equations, only spatially smooth noises
can be treated (see \cite{dBD1}, \cite{dBD2}).
Note that this equation is not strongly dissipative, indeed the eigenvalues of the linear
part do not grow to infinity.
However, it is known that Foias-Prodi type estimates
hold for the deterministic damped NLS equation (see \cite{goubet}) and we will see that
these can be generalized to the stochastic case and it is reasonable to think that
the above ideas can be generalized.

In this article, we show that the method based on
coupling argument is applicable. However it requires substantial adaptations. For instance,
contrary to the strongly dissipative
case treated in the above-mentioned articles, we are only able to prove a weaker form
of the Foias-Prodi estimates. Indeed, here, we prove that it holds in average and not path-wise. This
causes many technical difficulties when trying to use the coupling method. Moreover, another
important ingredient in the argument is an exponential estimate on the growth of the solution, which
we are unable to prove in our case. This is due to the fact that the Lyapunov structure is more
complicated here. It is not a quadratic functional.
We only prove polynomial estimate on the growth of the solutions and it results
that we can only prove
that convergence to equilibrium holds with polynomial speed at any order. Thus, we develop a
general result, which gives sufficient conditions for
polynomial mixing.

Note also that a crucial step in \cite{KS3} is the fact that the probability that  a solution enters
a ball of small radius is controlled precisely. This fact is still true for the damped NLS equation
considered here. However, its proof is more difficult than in the case of the Navier-Stokes equations
(see Proposition \ref{Prop_petit} and section 4 hereafter).

The remaining of  the article is divided into four parts. First, we give the notations,
and state our main result. Its proof is given in section 2. Section 3, 4 and 5 are devoted
to the proofs of intermediate results.



\section{Notation and Main result}
\

We set
$$
A = -\Delta, \quad  D(A) = H^1_0(0,1) \cap H^2(0,1)
$$
and write problem (\ref{EqIntro}) in the form
\begin{eqnarray}
\NLS{u} & = & b dW ,\label{Eqbase} \\
u(0) & = & u_0,  \label{Eqinitial}
\end{eqnarray}
where $W$ is a cylindrical Wiener process on $L^2(0,1)$ and $b$ is a linear operator on
$L^2(0,1)$.

We denote by  $(\mu_n)_n$ the increasing sequence of eigenvalues  of $A$ and
by $(e_n)_{n\in \N}$   the associated eigenvectors. Also, $P_N$ and $Q_N$
are the eigenprojectors onto the space $Sp(e_k )_{1\le k \leq n}$ and onto its complementary space.
Recall that for $s\ge 0$, $D(A^{s/2})$ is a closed subspace of $H^s(0,1)$ and that
$\|\cdot\|_s=|A^{s/2}\cdot|_{L^2(0,1)}$ is equivalent to the usual $H^s(0,1)$ norm on this space. Moreover
$$
D(A^{s/2})=\{u=\sum_{k\in\N}u_ke_k\in L^2(0,1)\;|\; \sum_{n\in\N} \mu_k^s u_k^2<\infty\}
\mbox{ and }\|u\|_s=\sum_{n\in\N} \mu_k^s u_k^2.
$$
We denote by $\abs{\cdot}$, $\abs{\cdot}_p$, $\norm{\cdot}$  the norms of
$L^2(0,1)$, $L^p(0,1)$, $H^1_0(0,1)$.

 The operator $b$ is supposed to commute with $A$, therefore it is  diagonal in the basis
 $(e_n)_{n\in\N}$ and we have
$$b e_n = b_n e_n.$$
We  assume that $b$ is Hilbert-Schmidt from $L^2(0,1)$ with values in
$D(A^{3/2})$. For any $s\in [0,3]$, we set
$$
B_s=\abs{b}^2_{\mathcal L_2(L^2(0,1),D(A^{s/2}))}=\sum_{n=0}^\infty \mu_n^s b_n^2.
$$

To study ergodic properties, we assume that there exists $N_*$ such that
\begin{equation}
\label{H}
 b_n>0, \mbox{ for  }n\leq N_*.
\end{equation}

The Hamiltonian plays an important role in the study of the nonlinear Schr\"odinger equation. It is
a conserved quantity in the absence of noise and damping. It is given by
$$
 \Hcal_*(v)=\frac{1}{2}\norm{v}^2-\frac{1}{4}\abssig{v}, \;  v\in H_0^1(0,1).
$$
In our study, it is the basic tool to derive a priori estimates. Recall that the Gagliardo-Nirenberg inequality gives a constant $c_0>0$ such that
$$
\abs{v}^4_4\leq \frac{1}{4}\norm{v}^2+\frac{c_0}{2}\abs{v}^6, \;  v\in H_0^1(0,1).
$$
It follows that, setting
$$
 \Hcal=\frac{1}{2}\norm{\cdot}^2-\frac{1}{4}\abssig{\cdot}+c_0\abs{\cdot}^6,
$$
we have
\begin{equation}
\label{e1.4}
\Hcal(v)\geq \frac{1}{4}\norm{v}^2+\frac{1}{4}\abssig{v } + \frac{c_0}{2}\abs{v}^6, \;  v\in H_0^1(0,1).
\end{equation}
In our computations, we will also use the following quantities which involve the $k^{\mbox{th}}$
power of the energy:
$$
E_{u,k}(t,s)=\Hcal(u(t))^k+\alpha k \int_s^t\Hcal(u(\sigma))^k d\sigma,\; t\ge s,
$$
when there is no ambiguity we set $E_{u,k}(t)=E_{u,k}(t,s)$.

In the following, $\alpha$, $B_s$ for $s\in [0,3]$ are fixed. All the constants appearing below may depend on them
 as well as on $A$, $b$.

Well-posedness of equations \eqref{Eqbase}, \eqref{Eqinitial} is easily proved.
Indeed, let $S(t)=e^{-\i At-\alpha t}$, $t\in \R$, be the group generated by the linear equation. We look for a mild solution,
that is a process $u$ with paths in $C(\R^+;H^1_0(0,1))$ satisfying
$$
u(t)= S(t) u_0+\i \int_0^tS(t-s)\abs{u(s)}^2u(s)ds+\int_0^t S(t-s) b dW(s).
$$
Since $(S(t))_{t\geq 0}$ is a contraction semi-group on $H^1_0(0,1)$ and the linear term is locally Lipschitz,
local in time existence and uniqueness is straightforward. Note that $\int_0^t S(t-s) b dW(s)$ lives in $D(A^\frac{3}{2})$.
An  a priori estimate is obtained thanks to Ito formula applied to $\Hcal$ and thanks to \eqref{e1.4}.
This use of Ito formula is not rigorous since $Au$ is not sufficiently smooth.
However, an approximation argument
can be used to prove rigorously this point.
For instance, the initial data can be approximated by a sequence in $D(A)$ and it is classical that if the initial
 data is in $D(A)$ then the solution is continuous with values in $D(A)$.

Note that in the following and especially in section $4$ and $5$, several computations are not rigorous due to the lack of
 regularity of the solutions. The same approximation argument should be applied.

By classical arguments, the solutions are strong Markov processes.
We denote by $(\Pcal_t)_{t\in \R^+}$ the
Markov transition semi-group associated to the solutions of \eqref{Eqbase}.

Also, given a Banach space $E$, the space $Lip_b(E)$ consists of all the bounded and Lipschitz
real valued functions on $E$. Its norm is given by
$$
\| \varphi \|_L= \| \varphi \|_\infty +L_\varphi, \; \varphi\in Lip_b(E),
$$
where $\|\cdot\|_\infty$ is the sup norm and $L_\varphi$ is the Lipschitz constant of $\varphi$.
The space of probability measures on $E$ is denoted by $\Pcal(E)$. It
can be  endowed with the metric defined by the total variation
 $$
\norm{\mu}_{var}= \sup  \left\{\abs{\mu(\Gamma)}\;|\; \Gamma \in \Br(E) \right\},
$$
where we denote by $\Br(E)$ the set of the Borelian subsets of $E$. It is well known that $\norm{.}_{var}$ is the dual
 norm of $\abs{.}_\infty$. We can also use a Wasserstein type metric
$$
\| \mu-\lambda\|_W=\sup_{\varphi\in Lip_b(E),\; \| \varphi\|_L\le 1}\abs{\int_E \varphi(u)d(\mu-\lambda)(u)}
$$
which is the dual norm of $\|\cdot\|_L$.
We also use the notation $\Dr(Z)$ for the distribution of a random variable $Z$.

The aim of this article is to establish the following result
\begin{Theorem}
\label{Th_MAIN}
There exists $N_0$ such that, if  (\ref{H}) holds with $N_*\geq N_0$, then
there exists a unique stationary
 probability measure $\nu$ of $(\Pcal_t)_{t\in \R^+}$ on $H_0^1(0,1)$.
 Moreover, for any  $p\in \N\backslash\{0\}$, $\nu$ satisfies
\begin{equation}\label{Eq_MAIN_a}
\int_{H^1_0(0,1)} \norm{u}^{2p} d\nu(u) < \infty,
\end{equation}
and there exists $C_p>0$  such that for any $\mu \in \Pcal(H^1_0(0,1))$
\begin{equation}\label{Eq_MAIN_b}
\|\Pcal^*_t\mu-\nu\|_W\leq C_p \left( 1+t \right)^{-p}\left(1+ \int_{H^1_0(0,1)} \norm{u}^2 d\mu(u) \right).
\end{equation}
\end{Theorem}
\begin{Remark}
Note that the existence of a stationary  measure
is a byproduct of the proof of the mixing property. It could be proved
directly by the standard argument involving the Krylov-Boboliubov
theorem. However, this would require more a priori estimates on the
solutions of the stochastic nonlinear Schr\"odinger equation.
\end{Remark}

\begin{Remark}
 In many articles and books, a family $\left(W_p(\cdot,\cdot)\right)_{p\in [1,\infty)}$  of Wasserstein type
 metrics is used.
 Actually, given a polish space $(E,d_E)$, 
these metrics are defined by
$$
W_p(\mu,\lambda)=\inf\left(\int_{E^2} d_E(x,y)^p\,\P(dx,dy)\right)^\frac 1p\quad \textrm{ for any }\quad p\in[1,\infty),
$$
where the infimum is taken over all probability measures $\P$ on $E^2$ whose marginal laws are $\lambda$ and $\mu$.

Let $(E,\norm\cdot_E)$ be a separable Banach space. We set $\,d_E(x,y)=\norm{x-y}_E\wedge 1$.
Then $(E,d_E)$ is a polish space.
Moreover, $W_1$ is equivalent to $\norm\cdot_W$. We have chosen
not to use the notation $W_1$ because it might lead to some confusion with the usual notation $W$
for a Wiener process.
\end{Remark}

The proof of our result is based on coupling arguments. These arguments have initially been used
in the context of stochastic partial differential equations in \cite{KS},
 \cite{KS2}, \cite{KS3},  \cite{Matt} . The main
difficulty here is that the nonlinear Schr\"odinger equation is not strongly dissipative and
several modifications are needed.

The strategy is the following. If the noise is non-degenerate, we observe that starting from
different initial data $u_0^1$, $u_0^2$, Girsanov transform can be used to
show that there exists a coupling $(u_1,u_2)$ of the law of the solutions
$u(\cdot,u_0^1)$, $u(\cdot, u_0^2)$ such that, for some time $T$, $u_1(T)=u_2(T)$
with positive probability. Iterating this argument, exponential convergence to equilibrium follows
(see section 1.1 in \cite{ODASSO}).
Here, as well as in the references above, the noise is assumed to be
 non-degenerate in the low modes only
$e_k$, $1\le k\le N_*$ and this argument gives a coupling such that $P_{N_*}u_1(T)=P_{N_*}u_2(T)$
with positive probability.
Another ingredient is used. It is based on the observation that if two solutions are such that their
low modes have been equal for a long time then they are very close
(see section 1.1 in \cite{ODASSO}).
In the case of a parabolic
equation, this is known as Foias-Prodi estimate. This can be generalized to dispersive equations
such as the Schr\"odinger equation considered here.
In \cite{goubet} this has been used to prove a property
of asymptotic smoothing in the deterministic case.

The main difference with the result in the
parabolic case is that we are not able to prove a path-wise Foais-Prodi estimate, we only
prove that this property holds in average. We need to introduce a substantial change in
the construction of the coupling. (See Remark \ref{changement}).
Moreover, here we only get polynomial convergence to equilibrium. This comes from the
fact that the Lyapunov functional adapted to the nonlinear Schr\"odinger equation is more
complicated, it is not a quadratic functional. We are not able to get exponential estimates on
the growth of the
solutions.

\section{Proof of Theorem \ref{Th_MAIN}}

We define $G$ by
$$
D(G)=D(A), \quad G v= \alpha v + \i A v,
$$
and set
$$
\Espace
\begin{array}{l}
X=P_{N_*} u,\, Y=Q_{N_*} u,\, \beta=P_{N_*} W,\, \eta= Q_{N_*} W,\\
 \sig_l=P_{N_*} b P_{N_*}, \, \sig_h=Q_{N_*} b Q_{N_*},\\
f(X,Y)=-\i P_{N_*}\left( \abs{X+Y}^{2}(X+Y) \right),\\
g(X,Y)=- \i Q_{N_*}\left( \abs{X+Y}^{2}(X+Y) \right).
\end{array}
$$
Then the nonlinear Schr\"odinger equation has the form
\begin{equation}\label{CGL_Eq_un_neuf}
\left \{
\Espace \begin{array}{lcl}
d X+ G X dt +f(X,Y)dt & = & \sig_l d\beta , \\
d Y+ G Y dt +g(X,Y)dt & = & \sig_h d\eta  , \\
\lefteqn{ X(0)=x_0 , \quad Y(0)=y_0.}
\end{array}
\right .
\end{equation}
Clearly \eqref{H} states that $\sig_l$ is invertible. We set
\begin{equation}\label{CGL_Eq_un_neuf_ter}
 \sig_0=\|\sig_l^{-1}\|_{{\mathcal L}(P_{N_*}L^2(0,1))}^{-1}>0.
\end{equation}
Given two initial data $u_0^i=(x_0^i,y_0^i)$, $i=1,2$,
we will construct a coupling
$(u_1,u_2)=\left((X_1,Y_1),(X_2,Y_2)\right)$ of the laws of the two solutions
$u(\cdot,u_0^i)=(X(\cdot,u_0^i),Y(\cdot,u_0^i))$, $i=1,2$, of \eqref{CGL_Eq_un_neuf}.
Recall that $(u_1,u_2)$ is a coupling of the laws of $u(\cdot,u_0^i)$, $i=1,2$, if the
distribution of $u_i$ is the distribution of $u(\cdot,u_0^i)$.

In fact we are going to construct a coupling $(V_1,V_2)=\left((u_1,W_1),(u_2,W_2)\right)$ of \linebreak
$\left\{\Dr \left((u(\cdot,u_0^i),W)\right)\right\}_{i=1,2}$ such that $X_i=P_{N_*} u_i$, $ \eta_i=Q_{N_*}W_i$ satisfy
good properties. More precisely, we want that $X_1=X_2$ and $\eta_1=\eta_2$ for as many trajectories as possible.
Clearly, we obtain a coupling of $\Dr(u(\cdot,u_0^1))$ and $\Dr(u(\cdot,u_0^2))$ .
Since the noise may degenerate in the equation for $Y$, we are not able to require that $u_1=u_2$.
The difference between $Y_1=Q_{N_*}u_1$ and $Y_2=Q_{N_*}u_2$ will be controlled thanks to a Foias-Prodi estimate.
Note that $W$ is a cylindrical process in $L^2(0,1)$ and does not live in $L^2(0,1)$. This is not a problem.
 Indeed, it is well-known that $W\in C(\R^+; H^{-1}(0,1))$ a.s. and we consider its distribution in this space.


We  define an integer valued random process $l_0$ which is particularly convenient
when deriving properties  of the coupling:
$$
l_0(k)=\min\left\{l \in \{0,...,k\} \,|\, P_{l,k} \;\textrm{ holds }\right\},
$$
 where $\min \emptyset =\infty$ and
$$
(P_{l,k})
\left \{
\Espace
\begin{array}{l}
X_1(t)=X_2(t), \quad \eta_1(t)=\eta_2(t), \quad \forall\; t \in [lT,kT],\\
\Hcal_l\leq d_0 ,\\
E_{u_i,4}(t,lT)\leq \kappa+1+d_0^4+d_0^8+B (t-lT), \quad \forall\; t \in [lT,kT],\; i=1,2,
\end{array}
\right.
$$
where we have set
$$
\Hcal_k=\Hcal(u_1(kT))+\Hcal(u_2(kT)).
$$
 We say that $(X_1,X_2)$ are coupled at $kT$ if $l_0(k)\leq k$, in other words
 if $l_0(k)\not = \infty$.

The coupling constructed below will be such that, for any $q\in\N\backslash\{0\}$, the following two
properties hold
\begin{equation}\label{CGL_Eq_a}
\left\{
\Espace
\begin{array}{l}
\exists \; d_0, \; \kappa,\; B, \; T_q>0  \; \textrm{ such that for any } l\leq k  \; , \;  T\geq T_q,  \\
\P\left(l_0(k+1)\not=l \; | \; l_0(k)=l \right) \leq  \frac{1}{2}\left(1+(k-l)T\right)^{-q}.
\end{array}
\right.
\end{equation}
This says that the probability that the trajectories decouple is small. Moreover, the longer they have
been coupled, the smaller this probability is.

The second property is that, for any $R_{0},d_0>0$,
\begin{equation}\label{CGL_Eq_b}
\left\{
\Espace
\begin{array}{l}
\exists \; T^*(R_0,d_0) >0  \textrm{ and } p_{-1}(d_0)>0 \textrm{ such that for any }T\geq T^*(R_0,d_0)\\
\P\left(l_0(k+1)=k+1 \; | \; l_0(k)=\infty ,\; \mathcal H_k \leq R_0 \right)
 \geq p_{-1}(d_0).
\end{array}
\right.
\end{equation}
In other words, inside a ball, the probability that two trajectories get coupled is bounded below.

The construction can be done by induction. At each step, we construct
a probability space $(\Omega_0,\F_0,\P_0)$ and a measurable couple of functions
 $(\omega_0,u_0^1,u_0^2)\to(V_i(\cdot,u_0^1,u_0^2))_{i=1,2}$   such that, for any
$(u_0^1,u_0^2)$, $(V_i(\cdot,u_0^1,u_0^2))_{i=1,2}$ is a coupling of
$\left\{\Dr\left((u(\cdot,u_0^i),W)\right)\right\}_{i=1,2}$ on $[0,T]$. Recall that the processes $(V_i)_{i=1,2}$
 live in the space
 $C(0,T;H^1_0(0,1))\times C(0,T;H^{-1}(0,1))$. We first set
$$
u_i(0)=u_0^i,\quad W_i(0)=0, \quad i=1,2.
$$
Assuming that we have built $(u_i,W_i)_{i=1,2}$ on $[0,kT]$, then we take $(V_i)_{i=1,2}$  as above independent
 of $(u_i,W_i)_{i=1,2}$ on $[0,kT]$ and set
$$
\left(u_i(kT+t),W_i(kT+t)\right) = V_i(t,u_1(kT),u_2(kT))
$$
for any $t \in [0,T]$.

The construction of $(V_i)_{i=1,2}$ depends on whether $l_0(k)\le k$ or $l_0(k)=\infty$. The two cases
are treated separately in sections 2.5. We first
 state and prove the Foias-Prodi
estimates and
give some a priori estimates. We then
recall some results  on coupling and give a general result implying polynomial
mixing. Sections 3, 4 and 5 are devoted to the proof of some results used in the course
of the proof.


\subsection{The Foias-Prodi estimates}

We define for any $(u_1,u_2,r)\in H^1_0(0,1)$
$$
J_*(u_1,u_2,r)=\frac{1}{2}\norm{r}^2
-\frac{1}{4}\int_{[0,1]}\left((\abs{u_1}^2+\abs{u_2}^2)\abs{r}^2+\left(\Re((u_1+u_2)\bar r)\right)^2\right)dx,
$$
where $\Re(z)$ is the real part of the complex number $z$, and
$$
J(u_1,u_2,r)=J_*(u_1,u_2,r)+c_1 \left(\sum_{i=1}^2\Hcal(u_i)\right)\abs{r}^2.
$$
We infer from the Sobolev Embedding $H^1(0,1)\subset L^\infty(0,1)$ that there exists $c>0$ such that
$$
\int_{[0,1]}\left((\abs{u_1}^2+\abs{u_2}^2)\abs{r}^2+\left(\Re((u_1+u_2)\bar r)\right)^2\right)dx
\le c (\|u_1\|^2 +\| u_2\|^2)|r|^2.
$$
Therefore, by \eqref{e1.4}, there exists $c_1>0$ such that
$$
J(u_1,u_2,r)\geq \frac{1}{4}\norm{r}^2.
$$
We set
$$
l(u_1,u_2)=1+\sum_{i=1}^2 \Hcal(u_i)^{4}.
$$
For $N\geq 1$, given $u_1,u_2$, two solutions of \eqref{Eqbase}, we define $J_{FP}^N=J_{FP}^N(u_1,u_2)$ by
$$
J_{FP}^N(t)=J(u_1(t),u_2(t),r(t))
\exp\left(2\alpha t-\frac{\Lambda}{\mu_{N+1}^{\frac{1}{8}}} \int_0^t  l(u_1,u_2)ds \right),
$$
where $r=u_1-u_2$.
The following result will be proved in section 5. It is  the Foias-Prodi estimates
adapted to the nonlinear Schr\"odinger equation. It states that two solutions having the same low
modes are close. The main difference with similar results in the parabolic case is that we are
not able to derive a path-wise estimate. Moreover, we introduce a slight generalization to allow
the perturbation of the Wiener process by a drift in the low modes. This generalization is essential in our
argument below.
\begin{Proposition}\label{Prop_Foias_Prodi_Hun}
For any $\kappa_0>0$,  there exists $\Lambda>0$ depending only on $\kappa_0$ 
 such that for any $N\in \N\backslash\{0\}$, we
 have the following property:

Let $W_1$, $W_2$ be two  cylindrical Wiener processes and
$h$ be an adapted process with continuous paths in $P_N L^2(0,1)$. Let $u_1$ be a solution in $C(0,T;H^1_0(0,1))$ of
$$
\Espace
\left\{
\begin{array}{rcl}
\NLS{u_1} & = & b dW_1 +hdt, \\
u_1(0) & = & u_0^1,
\end{array}
\right.
$$
and $u_2$ be the solution of \eqref{Eqbase}, \eqref{Eqinitial} for $u_0=u_0^2$ and $W=W_2$.
Let $\tau$  be a stopping time and assume that
 \begin{equation}\label{1.2.1}
 P_N u_1 = P_N u_2, \quad Q_N  W_1=Q_N W_2
 \mbox{  on  }  [0,\tau],
 \end{equation}
and
 \begin{equation}\label{1.2.2}
\norm{h(t)}^2\leq \kappa_0 l(u_1(t),u_2(t))^{3/4} \mbox{  on  }  [0,\tau],
 \end{equation}
then we have
 \begin{equation}\label{1.2.3}
\E \left(J_{FP}^N(u_1,u_2)(t\wedge\tau)\right) \leq  J(u_0^1,u_0^2,r_0),\; t>0,
 \end{equation}
 where $r_0=u_0^1-u_0^2$.
\end{Proposition}
We deduce a very useful Corollary.
\begin{Corollary}\label{Cor_Foias_Hun}
For any $B, \, d_0,\, \kappa_0 >0$, there exists $N_1(B,\kappa_0)$
and $C^*(d_0)$ such that, with the notations of
 Proposition \ref{Prop_Foias_Prodi_Hun}, if \eqref{1.2.1} and \eqref{1.2.2} hold with $N \geq N_1$, and
 for some $\rho>0$,
 \begin{equation}\label{1.2.4}
    E_{u_i,4}(t)\leq  \rho+1+d_0^{4}+d_0^{8}+B t \textrm{ on } [0,\tau], \mbox{ for } i=1,2,
 \end{equation}
then  for any $u_0^1,\, u_0^2$ such that $d_0\geq\sum_{i=1}^2\Hcal(u_0^i)$ and for any $a \in \R$,
$$
\P\left(\norm{r(T)} > C^* \left(d_0\right) \exp \left (a-\frac{\alpha}{4}T+ \rho  \right )\textrm{ and }
T\leq \tau  \right)\leq  \exp \left (-a-\frac{\alpha}{4}T \right ).
$$
Moreover, there exists $c>0$ such that
$$
C^*(d_0)\leq c d_0 e^{cd_0^8}.
$$
\end{Corollary}
Then, integrating \eqref{1.2.3} in Proposition \ref{Prop_Foias_Prodi_Hun} and  applying the inequality
$$
1+x\leq C_\delta e^{\delta x}\quad \textrm{for any } x\geq 0,
$$
 we obtain the following
result which, in Section 3, ensures that the Novikov condition holds and allows the use of the Girsanov Formula.
\begin{Lemma}
\label{Prop_Novikov_Hun}
For any $B, \, d_0,\, \kappa_0 >0$ and any $a\in \R$, there exists $N_2(B,\kappa_0, a)$
and $C^{**}(d_0,B)$ such that, with the notations of
 Proposition \ref{Prop_Foias_Prodi_Hun}, if \eqref{1.2.1} and \eqref{1.2.2} hold with $N \geq N_2$ and
   \eqref{1.2.4} holds for some $\rho>0$,
we obtain that for any $T$
$$
\Espace
\begin{array}{c}
\P\left(\int_{T}^\tau  l(u_1(s),u_2(s))\norm{r(s)}^2ds
> C^{**}(d_0,B) \exp \left (a-\frac{\alpha}{2}T+ \rho  \right )\textrm{ and }
 T\leq \tau  \right)\\
\leq \exp \left (-a-\frac{\alpha}{2}T  \right ).
\end{array}
$$
provided $d_0\geq\sum_{i=1}^2\Hcal(u_0^i)$ holds. Moreover, there exists $c>0$ such that
$$
C^{**}(d_0,B)\le C(B) d_0 e^{cd_0^8}.
$$
\end{Lemma}

We set
\begin{equation}
\label{N_0}
N_0=\max(N_1,N_2).
\end{equation}


\subsection{A priori estimates}

We first give an estimate proven in section 4
on the growth of the solutions of the stochastic nonlinear Schr\"odinger
equation.

\begin{Proposition}
\label{Prop_majoration_energie_sous_lineaire_Hun}
Assume that $u$ is a solution of
\eqref{Eqbase}, \eqref{Eqinitial}
associated with a Wiener process $W$.
Then, for any $(k,p) \in (\N\backslash\{0\})^2$, there exists $C_k'$ and $K_{k,p}$  depending only
on $k$ and $p$ 
such that for any $\rho>0$ and $0 \leq T <\infty$
$$
\Espace
\begin{array}{lcl}
\P \left (  \sup_{t \in [0,T]}  \left(   E_{u,k}(t) - C_k't  \right )
\geq  \Hcal(u_0)^{k}+\rho \left(\Hcal(u_0)^{2k}+T \right)  \right)& \leq & K_{k,p}\rho^{-p},\\
\P \left (  \sup_{t \in [T,\infty)} \left (   E_{u,k}(t) - C_k't  \right )
\geq   \Hcal(u_0)^{k}+\Hcal(u_0)^{2k}+1+\rho  \right) & \leq & K_{k,p}\left( \rho + T  \right)^{-p}.
\end{array}
$$
\end{Proposition}

The following result uses the Hamiltonian as a Lyapunov functional and is also proven in section 4.

\begin{Lemma}
\label{lem_Lyapounov_norm}
There exists $C_k>0$ such that for any $k\in\N\backslash\{0\}$, for any $t\in\R^+$ and for any stopping time $\tau$
$$
\Espace \left\{
\begin{array}{lcl}
\E\left( \mathcal H(u(t,u_0))^k \right)&\leq& \mathcal H(u_0)^k  e^{-\alpha k t } + \frac{C_k}{2},\\
\E\left( \mathcal H(u(\tau,u_0))^k 1_{\tau<\infty}\right)&\leq&  \mathcal H(u_0)^k +C_k\E(\tau 1_{\tau<\infty}).
\end{array}
\right.
$$
\end{Lemma}

The following result states that we control the probability of entering a small ball.

\begin{Proposition}
\label{Prop_petit}
 For any $R_0,R_1>0$, there exists $T_{-1}(R_0,R_1)\geq 0$ and $\pi_{-1}(R_1)>0$ such that
$$
\P\left( \Hcal(u(t,u_0^1))+\Hcal(u(t,u_0^2))\leq R_1 \right)\geq \pi_{-1}(R_1),
$$
provided $\Hcal(u_0^1)+\Hcal(u_0^2)\leq R_0$ and $t\geq T_{-1}(R_0,R_1)$.
\end{Proposition}


\subsection{Basic properties of couplings}
\

Let $(\mu_1,\mu_2)$ be two distributions on a same space $(E,\mathcal{E})$. Let $(\Omega,\mathcal{F},\P)$ be a probability
 space and let $(Z_1,Z_2)$ be two random variables $(\Omega,\mathcal{F}) \to (E,\mathcal{E})$. Recall
  that $(Z_1,Z_2)$
is a coupling of  $(\mu_1,\mu_2)$ if $\mu_i=\Dr(Z_i)$ for $i=1,2$ and that we have denoted by
$\Dr(Z_i)$ the law of the random variable $Z_i$.

 Let $\mu$, $\mu_1$ and $\mu_2$ be three probability measures on a space $(E,\mathcal{E})$ such that $\mu_1$ and $\mu_2$
are absolutely continuous
 with respect to $\mu$. We set
$$
d(\mu_1 \wedge \mu_2)=(\frac{d\mu_1}{d\mu}\wedge\frac{d\mu_2}{d\mu}) d\mu.
$$
This definition does not depend on the choice of $\mu$ and we have
$$
\norm{\mu_1-\mu_2}_{var}= \frac{1}{2} \int_E  \abs{ \frac{d\mu_1}{d\mu}-\frac{d\mu_2}{d\mu}}d\mu.
$$
Remark that if $\mu_1$ is absolutely continuous with respect to $\mu_2$, we have
\begin{equation}\label{e2.34bis}
\norm{\mu_1-\mu_2}_{var}\leq \frac{1}{2}\sqrt{\int \left ( \frac{d \mu_1}{d \mu_2} \right )^{2} d\mu_2-1}.
\end{equation}

Next result is a fundamental result in the coupling methods, the proof
is given for instance in the Appendix of \cite{ODASSO}.

\begin{Lemma}\label{lem_coupling}
Let $(\mu_1,\mu_2)$ be two probability measures on $(E,\mathcal{E})$. Then
$$
\norm{\mu_1-\mu_2}_{var}= \min \P(Z_1\not = Z_2).
$$
The minimum is taken over all couplings $(Z_1,Z_2)$ of $(\mu_1,\mu_2)$. There exists a coupling
which reaches the minimum value. It is called a maximal coupling
and has the following property:
$$
\P(Z_1=Z_2, Z_1 \in \Gamma)=(\mu_1 \wedge \mu_2)(\Gamma)\; \textrm{ for any } \Gamma \in \mathcal{E}.
$$
\end{Lemma}

\BLANC{
We also use the following result which is
lemma D.1 of \cite{Matt}.
\begin{Lemma}\label{lem_tech_coup_inf}
Let $\mu_1$ and $\mu_2$ be two probability measures on a space $(E,\mathcal{E})$. Let $A$ be an event of $E$. Assume that
 $\mu_1^A=\mu_1(A\cap.)$ is equivalent to $\mu_2^A=\mu_2(A\cap.)$.  Then for any $p>1$ and $C>1$
 $$
\int_A \left ( \frac{d \mu_1^A}{d \mu_2^A} \right )^{p+1} d\mu_2 \leq C <\infty
\quad \textrm{implies} \quad
\left(\mu_1\wedge\mu_2\right)(A) \geq  \left( 1-\frac{1}{p} \right)\left( \frac{\mu_1(A)^p}{pC} \right)^{\frac{1}{p-1}}.
$$
\end{Lemma}
}

Next result is a refinement of Lemma \ref{lem_coupling} used in \cite{Matt}
(see also Proposition 1.7  in \cite{ODASSO}).
\begin{Proposition}\label{Prop_Matt}
Let $E$ and $F$ be two Polish spaces, $f_0:E\to F$ be a measurable map  and $(\mu_1,\mu_2)$ be two probability measures
on $E$. We set
$$
\nu_i=f_0^* \mu_i, \quad i=1,2.
$$
Then there exist a coupling $(V_1,V_2)$ of $(\mu_1,\mu_2)$ such that $(f_0(V_1),f_0(V_2))$ is a maximal coupling of
 $(\nu_1,\nu_2)$.
\end{Proposition}



\subsection{Sufficient conditions for polynomial mixing}
\

We now state and prove a general
result which allows to reduce the proof of polynomial convergence to equilibrium
 to the verification of some conditions. This result is a polynomial version of Theorem 1.8 of subsection 1.3 in
 \cite{ODASSO} which gives
 sufficient conditions for exponential mixing.

We are concerned with
$v(\cdot,(u_0,W_0))=(u(\cdot,u_0),W(\cdot,W_0))$ a couple of
strongly Markov processes defined on Polish spaces $(E,d_E)$ and
$(F,d_F)$.
 We denote by $(\Pcal_t)_{t\in I}$ the Markov transition semigroup
 associated to $u$, where $I=\R^+$ or $T\N=\{kT, \; k\in \N\}$. We are also given a real
 valued function $\mathcal H$ defined on $E$.

We consider for any couple of initial conditions $(v_0^1,v_0^2)$  a coupling $(v_1,v_2)$ of
$\left\{\Dr(v(\cdot,v_0^1)),\Dr(v(\cdot,v_0^2))\right\}$. We write $v_i=(u_i,W_i)$.
Let  $l_0:\N\to \N\cup \{\infty\}$ be a random integer valued process which has
 the following properties
\begin{equation}\label{Abstract_ater}
\Espace
\left\{
\begin{array}{l}
l_0(k+1)=l \textrm{ implies } l_0(k)=l, \textrm{ for any } l\leq k,\\
l_0(k)\in \{0,1,2, ...,k\}\cup\{\infty\},\\
l_0(k) \textrm{ depends only of } v_1|_{[0,kT]} \textrm{ and } v_2|_{[0,kT]},\\
l_0(k)=k\; \textrm{ implies } \Hcal_k \leq d_0,
\end{array}
\right.
\end{equation}
where
$$
\mathcal H_k=\mathcal H(u_1(kT))+\mathcal H(u_2(kT)), \quad \mathcal H:\, E \to \R^+,
$$
and $d_0>0$.

We now give four conditions on the coupling. The first condition states that when $(v_1,v_2)$
 have been coupled for a long time then the probability that $(u_1,u_2)$ are close is high. It will be
 a consequence of the Foias-Prodi estimate.

\begin{equation}\label{Abstract_abis}
\Espace
\left\{
\begin{array}{l}
\textrm{There exist $c_0$ and $q>0$ such that  for any  $t \in [lT,kT]\cap I$}\\
\P\left(    d_E(u_1(t),u_2(t))> c_0 \left(t-lT\right)^{-q} \textrm{ and }l_0(k)\le l\right)
 \leq c_0 \left(t-lT\right)^{-q},
\end{array}
\right.
\end{equation}

The next two  properties are exactly \eqref{CGL_Eq_a} and \eqref{CGL_Eq_b}.

\begin{equation}\label{Abstract_a}
\left\{
\Espace
\begin{array}{l}
\exists  \;  (p_k)_{k \in \N},\; c_1>0,\; q_0>1+q  \; \textrm{ such that} , \\
\P\left(l_0(k+1)=l \; | \; l_0(k)=l \right) \geq p_{k-l}, \textrm{ for any } l\leq k , \\
1-p_k \leq  c_1 \left((k+1)T\right)^{-q_0},\; p_k>0  \;\textrm{ for any }  k \in \N.
\end{array}
\right.
\end{equation}
\begin{equation}\label{Abstract_b}
\left\{
\Espace
\begin{array}{l}
\textrm{There exist }  p_{-1}>0, \; R_0>0 \textrm{ such that}\\
\P\left(l_0(k+1)=k+1 \; | \; l_0(k)=\infty ,\; \mathcal H_k \leq R_0 \right)
 \geq p_{-1}.
\end{array}
\right.
\end{equation}

The last ingredient is the so-called Lyapunov structure and follows from Lemma \ref{lem_Lyapounov_norm}. It allows
the control of the probability to enter the ball of radius $R_0$. It states that there exist $K_1$ and $K'$ constants
such that for any initial data $v_0$
and any stopping time $\tau'$ taking values in $\{kT,\; k\in \N\}\cup\{\infty\}$
\begin{equation}\label{Abstract_c}
\Espace
\left \{
\begin{array}{lcl}
\E\mathcal H(v(t,v_0)) &\leq & e^{-\alpha t}\mathcal H(v_0)+\frac{K_1}{2},\quad t\geq 0, \\
\E\left(\mathcal H(v(\tau',v_0)) 1_{\tau'<\infty}\right) & \leq &
K' \left(\mathcal H(v_0) + \E(\tau' 1_{\tau'<\infty}) \right).
\end{array}
\right.
\end{equation}

The process $V=(v_1,v_2)$ is said to be $l_0$--Markov if the laws
of $V(kT+\cdot)$ and of $l_0(k+\cdot)-k$ on
 $\{l_0(k)\in \{k,\infty\}\}$
 conditioned by $\F_{kT}$ only depend on $V(kT)$ and are equal to the laws of $V(\cdot,V(kT))$ and $l_0$, respectively.

In this article, we construct a  coupling $(u_i,W_i)_{i=1,2}$ of
two solutions  which is $l_0$--Markov
 but not Markov. We could modify the construction so that it
 is Markov at discrete times $T\N=\{kT,\; k\in \N\}$. However,
 it does not seem to be possible to modify the coupling to be Markov at any times.
 As shown below, the following result implies Theorem \ref{Th_MAIN}. Its proof is given in section 3.

\begin{Theorem}\label{Th_Theo}
Assume that for any $(u_0^1,W_0^1)$, $(u_0^2,W_0^2)$ there exists
a coupling $V=(v_1,v_2)$ of the laws of
$(u(\cdot,u_0^1),W(\cdot,W_0^1))$ and
$(u(\cdot,u_0^2),W(\cdot,W_0^2))$ which is  $l_0$--Markov and
satisfies \eqref{Abstract_ater}, \eqref{Abstract_abis},
\eqref{Abstract_a}, \eqref{Abstract_b} and \eqref{Abstract_c} with
$R_0>4 K_1$ and $R_0\geq d_0$. Then there exists $c_4>0$  such
that, for any $\varphi \in Lip_b(E)$ and any $u_0^1,\; u_0^2 \in
E$,
\begin{equation}\label{Un_49_ter}
\abs{\E \varphi(u(t,u_0^1))-\E \varphi(u(t,u_0^2))}
\leq c_4 \left(1+t\right)^{-q}\|\varphi\|_{L} (1+\Hcal(u_0^1)
+\Hcal(u_0^2) ).
\end{equation}
\end{Theorem}
\begin{Corollary}\label{c2.5}
Under the same assumptions as Theorem \ref{Th_Theo},
there exists a unique stationary probability measure $\nu$ of $(\Pcal_t)_{t\in I}$ on $E$.
It satisfies,
\begin{equation}\label{Eq1.49bis}
\int_E \mathcal H(u)d\nu(u) \leq \frac{K_1}{2}.
\end{equation}
Moreover
for any $\mu \in \Pcal(E)$
\begin{equation}\label{Un_49_bis}
\|\Pcal^*_t\mu-\nu\|_W\leq 2c_4 \left(1+t\right)^{-q}\left(1+ \int_E \mathcal H(u)d\mu(u) \right).
\end{equation}
\end{Corollary}

To prove Theorem \ref{Th_Theo}, we first note that it is sufficient to prove that, for any initial data $u_0^1$ and
$u_0^2$, the coupling satisfies
\begin{equation}\label{Abs_a}
\P\left( d_E(u_1(t),u_2(t))>c_3 \left(1+t\right)^{-q} \right)\leq c_3  \left(1+t\right)^{-q}\left( 1+\Hcal(u_0^1)+\Hcal(u_0^2) \right)
\end{equation}
where, as above, $v_i=(u_i,W_i)$.
Indeed we
have, since $(u_1,u_2)$ is a coupling of $\left\{\Dr(u(\cdot,u_0^1)),\Dr(u(\cdot,u_0^2)\right\}$,
$$
\Espace
\begin{array}{l}
|\E \varphi(u(t,u_0^1))-\E \varphi(u(t,u_0^2))|
=|\E \varphi(u_1(t))-\E(\varphi(u_2(t)))|\\
 \le L_\varphi c_3 (1+t)^{-q}
+ 2 \| \varphi \|_\infty \P\left( d_E(u_1(t),u_2(t))>c_3 \left(1+t\right)^{-q} \right)\\
\le L_\varphi c_3 (1+t)^{-q} + 2 \| \varphi \|_\infty c_3  \left(1+t\right)^{-q}\left( 1+\Hcal(u_0^1)+\Hcal(u_0^2) \right)
\end{array}
$$
so that \eqref{Un_49_ter} follows. The existence and uniqueness of a stationary measure is then
straightforward. Moreover, \eqref{Un_49_bis} is an easy consequence of \eqref{Un_49_ter}
and \eqref{Eq1.49bis} follows from \eqref{Abstract_c}.

\subsection{Construction of the coupling}

We first state the following result.
\begin{Proposition}\label{Th_iso}
There exists a measurable map
$$
\Phi :  C((0,T);P_{N_*} H^1_0(0,1))\times C((0,T);Q_{N_*} H^{-1}(0,1)) \times  H_0^1(0,1)     \to  C((0,T);Q_{N_*} H^1_0(0,1)),
$$
such that for any $(u,W)$ solution of \eqref{Eqbase} and \eqref{Eqinitial}
$$
Y   =  \Phi( X ,\eta,u_0)\quad \textrm{ on } [0,T],\mbox{ where } X=P_{N_*} u,\; Y=Q_{N_*} u,\; \eta= Q_{N_*} W.
$$
Moreover $\Phi$ is a non-anticipative functions of $(X,\eta)$.
\end{Proposition}
To prove this result, we note that the equation
$$
y(t)=S(t)y_0-\int_0^tS(t-s)g(x(s),y(s))ds+\int_0^tS(t-s)dz(s),
$$
can be solved by a fixed point argument in $C(0,T;H^1_0(0,1))$ for any deterministic functions $x\in C(0,T;P_{N_*}H^1_0(0,1))$
 and $z\in C(0,T;D(A^\frac{3}{2}))$. The last term is defined thanks to an integration by part. Clearly
$y=\Psi(x,z,y_0)$ for a measurable function $\Psi$. Thus $Y=\Psi(X,\sig_h\eta,Q_{N_*} u_0)$. We set
$\Phi(x,\widetilde z,u_0)=\Psi(x,\sig_h \widetilde z,Q_{N_*} u_0)$ for $\widetilde z$ such that $\sig_h\widetilde z \in
C(0,T;D(A^\frac{3}{2}))$ and $0$ otherwise. It is clear that $\Phi$ is not anticipative.

As already explained, the coupling $(u_1,u_2)$ is constructed by induction and we start by constructing a coupling
for two solutions $u(\cdot,u_0^i),\; i=1,2$ on an interval $[0,T]$. In fact, we construct three different
couplings. At time $kT$, we choose between these depending on whether $l_0(k)=\infty$ and
$\Hcal(u_1(kT))+\Hcal(u_2(kT))\le R_0$
 (case a)
or $\l_0(k)\le k$ (case b). In this latter case, $P_{N_*}u_1(kT)=P_{N_*}u_2(kT)$. In the third case,
$l_0(k)=\infty$ and $\Hcal(u_0^1)+\Hcal(u_0^2)> R_0$, we choose the trivial coupling.


{\bf Case a:}  $l_0(k)=\infty$ and
$\Hcal(u_0^1)+\Hcal(u_0^2)\le R_0$. We construct a coupling such that \eqref{CGL_Eq_b} holds.

\smallskip

In this case, we consider  $u_0^1,\; u_0^2$ such that
$\Hcal(u_0^1)+\Hcal(u_0^2)\le R_0$. The construction of  the coupling is done in two steps. We set
$$
  \mu_i=\Dr \left((u(\cdot,u_0^i),W)\right), \; \textrm{ on } [0,T_1], \quad i= 1,2.
$$
\smallskip

{\underline{Step 1:}}

\smallskip

We first prove that, for any $d_0>0$, there exist  $T_1(d_0)>0$, $R_1=R_1(d_0)>0$ and  a coupling
 $(\widetilde V_i(\cdot,u_0^1,u_0^2))_{i=1,2}$ of
$(  \mu_1,  \mu_2)$
such that for any $(u_0^1,u_0^2)$  satisfying $\sum_{i=1}^2 \Hcal(u_0^i)\leq R_1$ we have
\begin{equation}\label{CGL_Eq_b_b}
\P\left( \widetilde X_1(T_1,u_0^1,u_0^2)= \widetilde X_2(T_1,u_0^1,u_0^2),
\;\sum_{i=1}^2\Hcal(\widetilde u_i(T_1,u_0^1,u_0^2))\leq d_0
 \right)\geq \frac{1}{2},
\end{equation}
where
$$
\widetilde V_i(\cdot,u_0^1,u_0^2)=\left( \widetilde  u_i(\cdot,u_0^1,u_0^2),
\widetilde W_i(\cdot,u_0^1,u_0^2) \right),
\;   \widetilde X_i(\cdot,u_0^1,u_0^2)=P_{N_*}\widetilde u_i(\cdot,u_0^1,u_0^2), \; i=1,2.
$$

To construct $\widetilde V_i$ such that \eqref{CGL_Eq_b_b} holds, we take $R_1,T_1>0$
and we set
$$
\Espace
\begin{array}{l}
E=C((0,T); H^1_0(0,1))\times C((0,T);H^{-1}(0,1)),\\
F=C((0,T); P_{N_*} H^1_0(0,1))\times C((0,T);Q_{N_*} H^{-1}(0,1)),\\
f_0\left(u,W\right)= (P_{N_*}u, Q_{N_*}W)=(X,\eta),\\
\hat \mu_1= \Dr\left((u(\cdot,u_0^1)+\frac{T_1-\cdot}{T_1}P_{N_*} (u_0^2-u_0^1),W)\right) \textrm{    on  } [0,T_1],\\
\nu_i=f_0^*\mu_i,\quad
\hat \nu_1=f_0^*\hat \mu_1.
\end{array}
$$
We apply Proposition \ref{Prop_Matt} to $(E,F,f_0,(\hat \mu_1,\mu_2))$ and obtain
$( \hat V_1(\cdot,u_0^1,u_0^2), \widetilde V_2(\cdot,u_0^1,u_0^2))$ a coupling of
$\left(\hat \mu_1,  \mu_2 \right)$. Moreover, setting
$$
(\widetilde X_2,\tilde \eta_2)=f_0(\widetilde V_2(\cdot,u_0^1,u_0^2)),\;
(\hat X_1,\eta_1)=f_0(\hat V_1(\cdot,u_0^1,u_0^2)),
$$
$((\widetilde X_2,\tilde \eta_2),(\hat X_1, \eta_1))$ is a maximal coupling of
$\left(\hat \nu_1,  \nu_2 \right)$.

Finally, we set
$$
\widetilde V_1= \left( \hat u_1 - \frac{T_1-\cdot}{T_1}P_{N_*} (u_0^2-u_0^1),W_1  \right) \textrm{ on } [0,T_1], \textrm{ where }
   \hat V_1= \left( \hat u_1 ,W_1  \right).
$$
We also write
$$
\beta_1=P_{N_*}W_1, \; \widetilde V_1= \left( \widetilde u_1,W_1  \right), \;
\; \widetilde V_2= \left( \widetilde u_2,W_2  \right).
$$
To prove \eqref{CGL_Eq_b_b} we first remark that since
$\hat u_1(T_1)=\widetilde u_1(T_1)$ and $\hat X_1= P_{N_*}\hat u_1$, $\widetilde X_i
=P_{N_*}\widetilde u_i$, then
\begin{equation}\label{S1}
\begin{array}{l}
\P\left( \widetilde X_1(T_1)=
 \widetilde  X_2(T_1) \textrm{  and  } \sum_{i=1}^2 \Hcal(\widetilde u_i(T_1))^6\leq
\kappa'(\rho,T_1,R_1)\right) \\
\geq \P\left(  \hat X_1=  \widetilde X_2 \textrm{ on } [0,T_1] \textrm{  and  } \sum_{i=1}^2
E_{\widetilde u_{i},6}(t) \leq \kappa'(\rho,t,R_1)
 \textrm{ on } [0,T_1] \right),
\end{array}
\end{equation}
where
$$
\kappa'(\rho,t,R_1)=2\left(R_1^6+C'_6 t+\rho(R_1^{12}+t)\right), \; t>0,
$$
$\rho$ to be chosen below.

Let us consider $ \bar X_1$ the unique solution of
\begin{equation}\label{CGL_Eq_b_d_un}
\Espace
\left\{
\begin{array}{l}
d \bar X_1+ G \bar X_1 dt-\delta(t)
+1_{t\leq\tau}f(\bar X_1-\hat \delta,\Phi(\bar X_1-\hat \delta,\eta_1,u_0^1))dt  =  \sig_l d\beta_1 , \\
\lefteqn{ \bar X_1(0)=x_0^2,}
\end{array}
\right.
\end{equation}
where $\delta(t)=\left(\frac{T_1-t}{T_1}-\frac{1}{T_1}G\right)P_{N_*} (u_0^2-u_0^1)$,
$\hat \delta(t)=\frac{T_1-t}{T_1}P_{N_*} (u_0^2-u_0^1)$,
 and $\tau=\tau_1\wedge\tau_2$ where
$$
\Espace
\left\{
\begin{array}{l}
 \tau_1= \inf\left\{t\in[0,T_1]\;|\;  E_{\bar X_1-\hat \delta+\Phi(\bar X_1-\hat \delta
,\eta_1,u_0^1),6}(t)  >  \kappa'(\rho,t,R_1)\right\},\\
 \tau_2= \inf\left\{t\in[0,T_1]\;|\;  E_{\bar X_1+\Phi(\bar X_1,\eta_1,u_0^2),6}(t)  >  \kappa'(\rho,t,R_1)\right\}.
\end{array}
\right.
$$
Clearly, $\bar X_1= \hat X_1=P_{N_*}\widetilde u_1+\hat\delta$ on $[0,\tau]$.
We denote by $\lambda_1$ the distribution of
$ (\bar X_1,\eta_1)$ under the probability $\P$.
We set $\widetilde{\beta}_1(t)=\beta_1(t)+\int_0^t d(s) dt $ where
\begin{equation}\label{CGL_Eq_b_d_j}
 d(t)=\delta(t)+
1_{t\leq\tau}\sig_l^{-1}
\Espace
\left(
\begin{array}{r}
f(\bar X_1(t)-\hat \delta(t),\Phi(\bar X_1-\hat \delta,\eta_1,u_0^1)(t))\\
-f(\bar X_1(t),\Phi(\bar X_1,\eta_1,u_0^2)(t))
\end{array}
\right).
\end{equation}
Then $\bar X_1$ is a solution of
\begin{equation}\label{CGL_Eq_b_d_i_deux}
\Espace
\left\{
\begin{array}{lcl}
d \bar X_1+ G \bar X_1 dt  +1_{t\leq\tau}f(\bar X_1,\Phi(\bar X_1,\eta_1,u_0^2))dt& = & \sig_l d\widetilde{\beta}_1 , \\
\lefteqn{ \bar X_1(0)=x_0^2.}
\end{array}
\right.
\end{equation}
It is not difficult to see that since $\sig_l$ is bounded below and by the definition of
$\tau$, the Novikov condition
is satisfied:
$$
\E\left(\exp\left(\int_0^T |d(t)|^2dt\right)\right)<\infty
$$
 and the Girsanov formula
can be applied. Then we set
$$
d\widetilde{\P}= \exp\left(\int_0^T d(s) dW(s)-\frac{1}{2} \int_0^T \abs{d(s)}^2 dt  \right)d\P
$$
 and deduce  that $\widetilde{\P}$ is a probability under which $(\widetilde{\beta}_1,\eta_1)$ is a
 cylindrical Wiener process. We denote by $\lambda_2$ the law of $ (\bar X_1,\eta_1)$ under $\widetilde{\P}$.

 We prove below that
\begin{equation}
\label{e2.38bis}
\Espace
\begin{array}{l}
\P\left(  \hat X_1(t)\not =  \widetilde X_2(t) \textrm{  or  } \sum_{i=1}^2
E_{\widetilde u_{i},6}(t) > \kappa'(\rho,t,R_1)
 \textrm{ for some } t<T_1 \right)\\
  \le \|\lambda_1-\lambda_2\|_{var} +
  \P\left(E_{\widetilde u_{1},6}(\tau) \ge  \frac{1}{2}\kappa'(\rho,\tau,R_1) \right)+
\P\left(E_{\widetilde u_{2},6}(\tau) \ge  \frac{1}{2}\kappa'(\rho,\tau,R_1)\right)
\end{array}
\end{equation}
We choose
$$
\rho=8K_{6,1}
$$
in the definition of $\kappa'(\rho,t,R_1)$ and deduce from Proposition \ref{Prop_majoration_energie_sous_lineaire_Hun} that
\begin{equation}\label{CGL_Eq_b_d_jbis2}
\P\left(E_{\widetilde u_{1},6}(\tau) \ge  \frac{1}{2}\kappa'(\rho,\tau,R_1) \right)+
\P\left(E_{\widetilde u_{2},6}(\tau) \ge  \frac{1}{2}\kappa'(\rho,\tau,R_1)\right)
\le \frac{1}{4}.
\end{equation}
Moreover using \eqref{e2.34bis}, we obtain
$$
\|\lambda_1-\lambda_2\|_{var} \le  \frac{1}{2}\sqrt{\E \exp\left( c \int_0^T \abs{d(s)}^2 dt \right) -1},
$$
and then, for $T_1,\; R_1$ sufficiently small,
$$
\|\lambda_1-\lambda_2\|_{var} \le  2(R_1\left(T_1+1\right)(1+R_1^2)+\kappa'(\rho,T_1,R_1)).
$$
We choose
$$
T_1=R_1,
$$
and deduce
\begin{equation}\label{CGL_Eq_b_d_jbis}
\|\lambda_1-\lambda_2\|_{var} \le  cR_1\left(1+R_1^{11}\right).
\end{equation}
Taking into account \eqref{S1}, \eqref{e2.38bis}, \eqref{CGL_Eq_b_d_jbis2} and \eqref{CGL_Eq_b_d_jbis},
 we can choose $R_1^0>0$ sufficiently small such that for any $R_1\le R_1^0$
\begin{equation}\label{CGL_Eq_b_ebis}
\P\left(  \widetilde X_1(T_1)=
\widetilde  X_2(T_1) \textrm{  and  } \sum_{i=1}^2 \Hcal(\widetilde u_i(T_1))^6\leq
\kappa'(\rho,R_1,R_1)\right)\geq \frac{1}{2}.
\end{equation}
Remark that there exists $R_1(d_0)\in (0,R_1^0)$ such that $R_1\leq R_1(d_0)$ implies
$$
\left\{ \sum_{i=1}^2 \Hcal(\widetilde u_i(T_1))^6\leq
\kappa'(\rho,R_1,R_1)\right\}\subset\left\{ \sum_{i=1}^2 \Hcal(\widetilde u_i(T_1))\leq
d_0\right\},
$$
so that \eqref{CGL_Eq_b_b} follows.

It remains to prove \eqref{e2.38bis}. We write
$$
\Espace
\begin{array}{l}
\P\left( \hat X_1(t)\not =  \widetilde X_2(t) \textrm{  or  } \sum_{i=1}^2
E_{\widetilde u_{i},6}(t) > \kappa'(\rho,t,R_1)
 \textrm{ for some } t<T_1 \right)\\
 = \P\left(\hat X_1|_{[0,\tau]}\not = \widetilde X_2|_{[0,\tau]} \textrm{  or  } \sum_{i=1}^2
E_{\widetilde u_{i},6}(\tau) = \kappa'(\rho,\tau,R_1)
  \right)\\
\le \P\left(  \hat X_1|_{[0,\tau]}\not = \widetilde X_2|_{[0,\tau]}\right) +
\P\left(E_{\widetilde u_{1},6}(\tau) \ge  \frac{1}{2}\kappa'(\rho,\tau,R_1) \right)\\
+\P\left(E_{\widetilde u_{2},6}(\tau) \ge  \frac{1}{2}\kappa'(\rho,\tau,R_1)
 \right).
\end{array}
$$
Let $\bar X_2$ is the solution of equation \eqref{CGL_Eq_b_d_i_deux} where $\beta_1$
is replaced by $\beta_2=P_{N_*}W_2$ then, with the probability $\P$,
$\bar X_2$ has the same law as $\bar X_1$ under the probability $\widetilde \P$ and
$$
\P\left( P_{N_*} \hat u_1|_{[0,\tau]}\not = P_{N_*} \widetilde u_2|_{[0,\tau]}\right)
\le \P(\bar X_1\not = \bar X_2).
$$
Thus, \eqref{e2.38bis} would follow if $((\bar X_1,\eta_1),(\bar X_2,\eta_2))$ was a maximal coupling
of $(\lambda_1,\lambda_2)$ (here, we have set $\eta_2=Q_{N_*}W_2$). However, we only know
that $((\hat X_1, \eta_1),(\widetilde X_2,\tilde \eta_2))$ is a maximal coupling of
$\left(\hat \nu_1,  \nu_2 \right)$. It is not difficult to remedy this problem. Indeed, the above result
holds for any  coupling of $\left(\hat \nu_1,  \nu_2 \right)$. Thus, instead of
$((\hat X_1, \eta_1),(\widetilde X_2,\tilde \eta_2))$, we choose another coupling such that
the processes constructed as
$((\bar X_1,\eta_1),(\bar X_2,\eta_2))$ above is a maximal coupling
of $(\lambda_1,\lambda_2)$. Then, the right hand side is equal to the right hand side
of \eqref{e2.38bis} while, by Lemma \ref{lem_coupling}, the left hand side is larger than the left
hand side of \eqref{e2.38bis}.

\smallskip

{\underline{Step 2:}} Construction of the coupling under the assumptions of case a.

\smallskip

Thanks to Proposition \ref{Prop_petit}, we know that there exists
$T_{-1}(R_0,R_1)>0$ and $\pi_{-1}(R_1)>0$ such that
\begin{equation}\label{Eq_b_c}
\P\left(\sum_{i=1}^2 \Hcal(u(\theta,u_0^i))\leq R_1\right)\geq \pi_{-1}(R_1),
\end{equation}
provided $\sum_{i=1}^2 \Hcal(u_0^i)\le R_0$ and  $\theta\geq T_{-1}(R_0,R_1)$.

We set $T^*(R_0,d_0)=T_{-1}(R_0,R_1(d_0))+T_{1}(d_0)$ and assume that $T\geq T^*(R_0,d_0)$. We also write
$\theta=T-T_1$.
Then on $[0,\theta]$, we take the trivial coupling which we denote by $(V'_1,V'_2)$.
Finally, we consider
 $(\widetilde V_1,\widetilde V_2)$ as above independent of $(V'_1,V'_2)$ and we set
$$
\Espace
V^a_i(t,u_0^1,u_0^2)=\left\{
\begin{array}{ll}
V'_i(t,u_0^1,u_0^2) & \textrm{ if } t\leq \theta, \\
\widetilde V_i(t-\theta,V'_1(\theta,u_0^1,u_0^2),V'_2(\theta,u_0^1,u_0^2)) & \textrm{ if } t\geq \theta.
\end{array}
\right.
$$
Combining \eqref{CGL_Eq_b_b} and \eqref{Eq_b_c} and setting
$$
p_{-1}(d_0)=\frac{1}{2}\pi_{-1}(R_1(d_0)),
$$
we obtain, for any $u_0^1,\; u_0^2$ such that
$\Hcal(u_0^1)+\Hcal(u_0^2)\le R_0$,
\begin{equation}\label{CGL_Eq_b_c}
\P\left(  X_1^a(T,u_0^1,u_0^2)=  X^a_2(T,u_0^1,u_0^2),
\;\sum_{i=1}^2\Hcal( u^a_i(T,u_0^1,u_0^2))\leq d_0
 \right)\geq  p_{-1}(d_0),
\end{equation}
where now
$$
V^a_i(\cdot,u_0^1,u_0^2)=\left(   u_i^a(\cdot,u_0^1,u_0^2),   W_i^a(\cdot,u_0^1,u_0^2) \right),
\;    X_i^a(\cdot,u_0^1,u_0^2)=P_{N_*}u_i^a(\cdot,u_0^1,u_0^2), \; i=1,2.
$$
Clearly, \eqref{CGL_Eq_b_c} implies \eqref{CGL_Eq_b}.

\smallskip

{\bf Case b:} $l_0(k)\le k$. We now construct a coupling so that  \eqref{CGL_Eq_a} holds. Since
 \eqref{CGL_Eq_a} depends on the whole history of the coupling and not only on the latest step,
  \eqref{CGL_Eq_a} is proved afterwards when the coupling is constructed on $[0,\infty)$.

\smallskip

In this case, we have $P_{N_*}u_0^1=P_{N_*}u_0^2$.
We write $x=P_{N_*}u_0^1=P_{N_*}u_0^2$, $y_1=Q_{N_*}u_0^1$ and $y_2=Q_{N_*}u_0^2$.

We apply Proposition \ref{Prop_Matt} to
$$
\Espace
\begin{array}{l}
E=C((0,T); H^1_0(0,1))\times C((0,T);H^{-1}(0,1)),\\
F=C((0,T);P_{N_*} H^1_0(0,1))\times C((0,T);Q_{N_*} H^{-1}(0,1)),\\
f_0\left(u,W\right)=(P_{N_*}u,Q_{N_*}W)= (X,\eta), \\
\mu_i=\Dr\left((u(\cdot,u_0^i),W)\right), \quad \textrm{ on } [0,T].
\end{array}
$$
We set $\nu_i=f_0^*\mu_i=\Dr\left((X(\cdot,u_0^i),\eta)\right)$ on $[0,T]$.
We obtain $(V_i^b(\cdot,u_0^1,u_0^2))_{i=1,2}=(u_i^b(\cdot,u_0^1,u_0^2),
W_i^b(\cdot,u_0^1,u_0^2))_{i=1,2}$, a coupling of $(\mu_1,\mu_2)$ such that if
we set
$$
( X^b_i,\eta^b_i)=f_0(V^b_i), \quad i=1,2.
$$
Then $( X^b_i, \eta^b_i)(\cdot,u_0^1,u_0^2))_{i=1,2}$
is a maximal coupling of $(\nu_1,\nu_2)$. We define $Y^b_i=Q_{N_*}u^b_i$,
$\beta^b_i=P_{N_*}W_i^b$

Let $\tau$ be  a stopping time associated to the process $(X,\eta)$.

\BLANC{We wish to estimate $\norm{}_{var}$. This enables us to use Lemma \ref{lem_coupling}
and to obtain an estimate which will be crucial to prove \eqref{CGL_Eq_a}.}

Let    $\widetilde X_1^b$ be
the unique solution of the truncated equation
\begin{equation}\label{CGL_Eq_i_un}
\Espace
\left\{
\begin{array}{lcl}
d \widetilde X_1^b+ G \widetilde X_1^b dt +
1_{t\leq\tau}f(\widetilde X_1^b,\Phi(\widetilde X_1^b, \eta_1^b,(x,y_1)))dt  =  \sig_l d\beta_1^b , \\
\lefteqn{ \widetilde X_1^b(0)=x.}
\end{array}
\right.
\end{equation}
Clearly $\widetilde X_1^b= X_1^b$ on $[0,\tau]$.  We denote by $\lambda_1$  the distribution of
$(\widetilde X_1^b,\eta)$ under the probability $\P$.

Let
$\widetilde{\beta_1^b}(t)=\beta_1^b(t)+\int_0^t d(s) dt $ where
$$
 d(t)=1_{t\leq\tau}\left(\sig_l\right)^{-1}\left(f(\widetilde X_1^b(t),
 \Phi(\widetilde X_1^b, \eta_1^b,(x,y_2))(t))-
 f(\widetilde X_1^b(t), \Phi(\widetilde X_1^b, \eta_1^b,(x,y_1))(t))\right).
$$
We take below a stopping time  $\tau$  such that
\begin{equation}
\label{novikov}
\int_0^T|d(t)|^2dt\leq M,
\end{equation}
for a constant $M$. Thus Novikov condition holds and  Girsanov formula applies. Setting
$$
d\widetilde{\P}= \exp\left(\int_0^T d(s) dW(s)-\frac{1}{2} \int_0^T \abs{d(s)}^2 dt  \right)d\P,
$$
we know that $\widetilde{\P}$ is a probability under which $(\widetilde{\beta},\eta)$ is a
cylindrical Wiener process.

Furthermore, with such a stopping time  $\tau$,   $\widetilde X_1^b$ is the solution of
\begin{equation}\label{CGL_Eq_i_deux}
\Espace
\left\{
\begin{array}{l}
d \widetilde X_1^b+ G \widetilde X_1^b dt +1_{t\leq\tau}f(\widetilde X_1^b,
\widetilde \Phi(\widetilde X_1^b, \eta_1^b,(x,y_2)))dt  =  \sig_l d\widetilde{\beta} , \\
\lefteqn{ \widetilde X_1^b(0)=x.}
\end{array}
\right.
\end{equation}
We denote by $\lambda_2$ the law of $(\widetilde X_1^b,\eta)$ under $\widetilde{\P}$.
 As in the case a, it is not difficult to see that
\begin{equation}\label{CGL_Eq_jter}
\norm{\lambda_1-\lambda_2}_{var} \leq
\frac{1}{2}\sqrt{\E \exp\left( c \int_0^T \abs{d(s)}^2 dt \right) -1}.
\end{equation}
This will be helpful to estimate $\norm{\nu_1-\nu_2}_{var}$.

\smallskip

{\bf Definition of the coupling on $[0,\infty)$.}

\smallskip

We first set
$$
u_i(0)=u_0^i,\quad W_i(0)=0, \quad i=1,2.
$$
Assuming that we have built $(u_i,W_i)_{i=1,2}$ on $[0,kT]$, then we take $(V^a_i)_i$ and
$(V^b_i)_i$ as above independent
 of $(u_i,W_i)_{i=1,2}$ on $[0,kT]$ and set for any $t \in [0,T]$
$$
\left(u_i(kT+t),W_i(kT+t)\right)=
\quad\quad\quad\quad\quad\quad\quad\quad\quad\quad\quad\quad\quad\quad\quad\quad\quad
\quad\quad\quad\quad\quad
$$
\begin{equation}\label{CGL_Eq_m_bis}
\Espace
\left \{
\Espace \begin{array}{ll}
V^a_i(t,u_1(kT),u_2(kT))& \textrm{ if } l_0(k) =\infty \textrm{ and } \Hcal(u_0^1)+\Hcal(u_0^2)\le R_0,\\
V^b_i(t,u_1(kT),u_2(kT)) & \textrm{ if } l_0(k) \leq k,\\
V^0_i(t,u_1(kT),u_2(kT))& \textrm{ if } l_0(k) =\infty \textrm{ and } \Hcal(u_0^1)+\Hcal(u_0^2)> R_0,
\end{array}
\right.
\end{equation}
where $V^0_i(t,u_1(kT),u_2(kT))$ is the trivial coupling. In other words, we take a cylindrical Wiener process $W$
 independent of $(u_i,W_i)_{i=1,2}$ on $(0,kT)$ and set
$$
V^0_i(t,u_1(kT),u_2(kT))=\left((u(t-kT,u_0^1),W),(u(t-kT,u_0^2),W)\right).
$$
Remark that, when $l_0(k) =\infty \textrm{ and } \Hcal(u_0^1)+\Hcal(u_0^2)> R_0$, the choice of the coupling is not very important.

 Clearly,
$(u_i,W_i)_{i=1,2}$ is a coupling of $(u(\cdot,u_0^i))_{i=1,2}$
which is $l_0$--Markov. In the following, we write
$$
X_i=P_{N_*}u_i,\; Y_i=Q_{N_*}u_i,\; \beta_i=P_{N_*}W_i,\;  \eta_i=Q_{N_*}W_i,\; i=1,2.
$$

It remains to prove that  \eqref{CGL_Eq_a} holds.

{\bf Proof of \eqref{CGL_Eq_a}}

\smallskip

We are in the situation where the coupling on $[kT,(k+1)T]$ has been constructed in
case b. We use the notation used in the construction of the coupling.

 Let us define for $i=1,2$
$$
\hat \tau^i_{k,l}=\inf\left\{t\in[0,T]\;\left| \;
 E_{\hat u_{i,4}}(kT+t,lT)>\kappa+1+d_0^4+d_0^8+C'_4 (t+(k-l)T)\right.\right\},
 $$
where $C'_4$ is given in Proposition \ref{Prop_Foias_Prodi_Hun}, and
 $$
 \hat\tau^3_{k,l}=\inf\left\{ t\le  T\;\left| \;
 \int_{kT}^{kT+t\wedge\hat\tau^1_{k,l}\wedge  \hat\tau^2_{k,l}}
  l(\hat u_1(s),\hat u_2(s))\norm{\hat r(s)}^2 ds > C_*(d_0) e^{a-\frac{\alpha}{2} (k-l)T}
\right.
\right\},
$$
where $a,\; d_0, \;\kappa$ are chosen below, $C_*(d_0) =C^{**}(C'_4,d_0)$ is given in Lemma
\ref{Prop_Novikov_Hun} and
$$
\Espace
\begin{array}{l}
 \hat u_i = u_i \; \textrm{  on } [0,kT], \quad \hat u_i(kT+\cdot)
 =\widetilde X_1^b+\Phi(\widetilde X_1^b,\eta_1^b,u_i(kT))\; \textrm{  on } [kT,(k+1)T],\\
\hat r=\hat u_1-\hat u_2.
\end{array}
$$
We also take $B=C'_4$ in the definition of $l_0(k)$.

Note that, with the notation of case b, $\hat u_1$ (resp. $\hat u_2$) is  a solution of a truncated
NLS equation under the the probability $\P$ (resp. $\widetilde \P$). It follows that when
$(\widetilde X_1^b, \eta_1^b)$ has law $\lambda_1$ (resp. $\lambda_2$) then $\hat u_1$
(resp. $\hat u_2$)
is a solution in law of a truncated NLS equation. But if $(\widetilde X_1^b, \eta_1^b)$ has law
$\lambda_1$,  $\hat u_2$ is a solution of a truncated NLS equation with a drift term.

We wish to use the construction described in case b with the stopping time $\tau=\tau_{k,l}$ given by
$$
\tau_{k,l}= \hat \tau^1_{k,l}\wedge \hat \tau^2_{k,l}\wedge \hat \tau^3_{k,l}.
$$
Then
$$
| d(t)|\le 1_{t\leq\tau_{k,l}}\sigma_0|f(\widetilde X_1^b(t),
 \Phi(\widetilde X_1^b, \eta_1^b,(x,y_2))(t))-
 f(\widetilde X_1^b(t), \Phi(\widetilde X_1^b, \eta_1^b,(x,y_1))(t))|
$$
and it is not difficult to see that
$$
| d(t)|^2\le c1_{t\leq\tau_{k,l}}  l(\hat u_1(t),\hat u_2(t))\|\hat r(t)\|^2.
$$
So that, by the definition of $\tau_{k,l}$, we get
\begin{equation}\label{CGL_Eq_i_trois}
\int_0^T\abs{d(t)}^2 dt\leq C^*(d_0) \sig_0^{-2}  \exp\left( a-\frac{\alpha}{2}(k-l)T \right).
\end{equation}
Hence the Novikov condition is satisfied and \eqref{e2.34bis} holds.

Moreover, using the same argument as in the proof of \eqref{e2.38bis}, we obtain
\begin{equation}
\label{e2.38}
\P\left((X_1^b,\eta_1^b)\not =(X_2^b,\eta_2^b) \textrm{ or } \tau<T\right)\leq
 \|\lambda_1-\lambda_2\|_{var} +\nu_1(\hat A^c_1)+ \nu_1(\hat A^c_3) + \nu_2(\hat A^c_2).
\end{equation}
where
$$
 \hat A_i=\{(X,\eta )\;|\;\hat \tau_i=T\},\; i=1,2,3.
$$
It can be seen that for $ i=1,2$
\begin{equation}
\label{e1}
\nu_i(\hat A^c_i)=\P\left( \sup_{t\in [0,T]}\left(E_{u_i,4}(kT+t,lT)-C'_4(t+(k-l)T)\right)>\kappa+1+d_0^4
+d_0^8 \, | \, \F_{kT}\right).
\end{equation}
The third term $\nu_1(\hat A^c_3)$ cannot be written in terms of $u_1$ and $u_2$. Indeed,
when $(\widetilde X_1^b, \eta_1^b)$ has law $\nu_1$, $\hat u_2$ is a solution of an equation
with a drift term.
\begin{Remark}
\label{changement} We remark here that Proposition
\ref{Prop_Foias_Prodi_Hun} is not the  Foias-Prodi estimate which
is usually used in the coupling method. Here, we have also a drift
term $h$. This modification is introduced precisely to treat the
term $\nu_1(\hat A^c_3)$. We take $h(\cdot)=b d(kT+\cdot)=\sig_l
d(kT+\cdot)$. This additional term is due to the fact that we
introduce a term depending on $r$ in the truncation. In the
preceding papers using this kind of coupling method, this was not
necessary and the Foias-Prodi estimate was used to get
\eqref{CGL_Eq_i_trois}. However, this requires a path-wise
Foias-Prodi estimate and we do not know if it holds in our
situation.
\end{Remark}
By \eqref{e1}, we have
$$
\nu_1(\hat A^c_1)+ \nu_1(\hat A^c_3) + \nu_2(\hat A^c_2)\le 3\P(B_{l,k}\, \big|\, \F_{kT})
$$
with
$$
B_{l,k}=\left\{
\begin{array}{l}
\sup_{t\in [0,T]}\left(E_{u_i,4}(kT+t,lT)-C'_4(t+(k-l)T)\right)>\kappa+1+d_0^4
+d_0^8,\;i\in\{1,2\}  \\
\textrm{ or}\int_{kT}^{kT+\tau_{k,l}} \left(\sum_{i=1}^2\Hcal( \hat u_i(s))^2 \right)\norm{\hat r(s)}^2 ds \ge C_*(d_0) e^{2\kappa-\frac{\alpha}{4} (k-l)T}
\end{array}
\right\}.
$$
Let us write
$$
\begin{array}{l}
\P(B_{l,k}\,\big|\,l_0(l)=l)\\
\le \ds{\sum_{i=1,2}
\P\left(\sup_{t\in [0,T]}\left(E_{u_i,4}(kT+t,lT)-C'_4(t+(k-l)T)\right)>\kappa+1+d_0^4+d_0^8\,\big|\,l_0(l)=l\right)}\\
\ds{+ \P\left(\int_{kT}^{kT+\tau_{k,l}} \left(\sum_{i=1}^2\Hcal( \hat u_i(s))^2 \right)\norm{\hat r(s)}^2 ds \ge C_*(d_0) e^{2\kappa-\frac{\alpha}{4} (k-l)T}\,\big|\,l_0(l)=l\right)}.
\end{array}
$$
Using Proposition \ref{Prop_majoration_energie_sous_lineaire_Hun} with $\kappa=\rho$
and solutions starting at
$lT$ and replacing $T$ by $kT$ we get, since $l_0(l)=l$ implies $\Hcal(u_i(lT)\le d_0$,
$$
\begin{array}{l}
\ds{ \P\left(\sup_{t\in [0,T]}\left(E_{u_i,4}(kT+t,lT)-C'_4(t+(k-l)T)\right)>\kappa+1+d_0^4+d_0^8 \, \big|\, \F_{lT}\right)}\\
\le \ds{ \P\left(\sup_{t\in [0,T]}\left(E_{u_i,4}(kT+t,lT)-C'_4(t+(k-l)T)\right)>\kappa+1+\Hcal(u_i(lT)^4+
\Hcal(u_i(lT)^8 \, \big|\, \F_{lT}\right)}\\
\ds{\le K_{4,q+1} \left( \kappa+(k-l)T) \right)^{-q-1}}.
\end{array}
$$
It follows
$$
\begin{array}{l}
\ds{ \P\left(\sup_{t\in [0,T]}\left(E_{u_i,4}(kT+t,lT)-C'_4(t+(k-l)T)\right)>\kappa+1+d_0^4+d_0^8 \, \big|\, l_0(l)=l\right)}\\
\ds{\le K_{4,q+1} \left( \kappa+(k-l)T) \right)^{-q-1}.}
\end{array}
$$
Similarly, by Lemma \ref{Prop_Novikov_Hun}, with $h(t)=\sig_l d(kT+t)1_{t\le \tau}$
which clearly satisfies \eqref{1.2.2} and $\rho=a=\kappa$, we have
$$
\begin{array}{l}
\ds{ \P\left(\int_{kT}^{kT+\tau_{k,l}} \left(\sum_{i=1}^2\Hcal( \hat u_i(s))^2 \right)\norm{\hat r(s)}^2 ds \ge C_*(d_0) e^{2\kappa-\frac{\alpha}{4} (k-l)T}\,\big|\,l_0(l)=l\right)}\\
\\
\ds{\le e^{-\kappa -\frac\alpha2(k-l)T}}\\
\le c\left( \kappa+(k-l)T) \right)^{-q-1}.
\end{array}
$$
Gathering these estimates, we obtain
$$
\P(B_{l,k}\,\big|\,l_0(l)=l)\le c\left( \kappa+(k-l)T) \right)^{-q-1}.
$$
By \eqref{e2.38}, \eqref{e2.34bis}, and \eqref{CGL_Eq_i_trois}, we obtain
for  $k\ge l$ and on $l_0(k)=l$
$$
\Espace
\begin{array}{l}
\P\left( (X_1,\eta_1) \not = (X_2,\eta_2) \textrm{ on } [kT,(k+1)T] \textrm{ or } B_{k,l} \,\left| \, \F_{kT} \right.\right)\\
\le \|\lambda_1-\lambda_2\|_{var}+ 3\P(B_{l,k} \big|  \F_{kT} )\\
\le \sqrt{\E \exp\left( c \int_0^T \abs{d(s)}^2 dt \right) -1} +3\P(B_{l,k}\big|\F_{kT})\\
\le  C^*(d_0) \sig_0^{-1}  e^{ \kappa-\frac{\alpha}{4}(k-l)T }
+3\P(B_{l,k}\big|\F_{kT}).
\end{array}
$$
We have
$$
\{l_0(k)=l\}\cap\{(X_1,\eta_1)  = (X_2,\eta_2) \textrm{ on } [kT,(k+1)T]\}\cap B_{l,k}^c\subset \{l_0(k+1)=l\}.
$$
Therefore, integrating over $l_0(k)=l$ gives for $T\geq T_1(d_0)$ and for $k>l$
$$
\P\left( l_0(k+1)\not=l,\; l_0(k)=l \,|\,l_0(l)=l\right)
\leq C^*(d_0) \sig_0^{-1}  e^{ \kappa-\frac{\alpha}{4}(k-l)T }
+3\P(B_{l,k}\,|\,l_0(l)=l).
$$
Which implies that there exists $\kappa>0$ sufficiently large and $d_0>0$ sufficiently small such that for any $T>0$
\begin{equation}
\label{e2.39}
\P\left( l_0(k+1)\not=l ,\; l_0(k)=l\,|\,l_0(l)=l\right)\leq \frac{1}{4}\left( 1+(k-l)T) \right)^{-q}.
\end{equation}
Remark that
$$
\P\left( l_0(k)\not=l|l_0(l)=l \right) \leq \sum_{n=l}^{k-1} \P\left( l_0(n+1)\not=l,\; l_0(n)=l    \,|\,l_0(l)=l\right),
$$
so that, applying \eqref{e2.39}, we obtain
$$
\P\left( l_0(k)\not=l|l_0(l)=l \right) \leq \frac{1}{4}+ \frac{1}{T^{q+1/2}}
\sum_{n=1}^\infty \frac{1}{k^{q}}\leq \frac{1}{4} +
C_q \frac{1}{T^{q}},
$$
which implies that there exists $T_q>0$ such that for $T\geq T_q$
\begin{equation}\label{CGL_Eq_eter}
\P\left( l_0(k)=l|l_0(l)=l \right) \geq \frac{1}{2},
\end{equation}
 Combining \eqref{e2.39} and \eqref{CGL_Eq_eter}, we establish \eqref{CGL_Eq_a}.


\subsection{Conclusion}

We have just shown that the coupling constructed in section 2.5 satisfies \eqref{CGL_Eq_a}
and \eqref{CGL_Eq_b} which are precisely \eqref{Abstract_a} and \eqref{Abstract_b}. The constants
used in \eqref{CGL_Eq_a} have been chosen in the preceeding subsection. The random
variables $l_0(k)$ clearly satisfy \eqref{Abstract_ater} and, as already mentioned, \eqref{Abstract_c}
is implied by Lemma \ref{lem_Lyapounov_norm}. Finally, \eqref{Abstract_abis} is a consequence
of Proposition \ref{Prop_Foias_Prodi_Hun} with $h=0$ and Tchebychev inequality.

We deduce that Theorem \ref{Th_Theo} can be applied. Moreover \eqref{Eq_MAIN_a}
is a consequence of Lemma
\ref{lem_Lyapounov_norm}.
This ends the proof of Theorem \ref{Th_MAIN}.


\section{Proof of Theorem \ref{Th_Theo}}
\

\subsection{Reformulation of the problem}

We already noticed that
it is sufficient to establish \eqref{Abs_a}.

 Let us  denote by $k$ the unique integer such that \mbox{$t\in (2(k-1)T,2kT]$}. Then
$$
\Espace
\begin{array}{l}
\P\left( d_E(u_1(t),u_2(t))>c_0  \left(1+t-(k-1)T\right)^{-q}\right)\leq \P\left(l_0(2k)\geq k\right)+
\\
\quad \quad \quad \quad \quad \quad \quad \P\left( d_E(u_1(t),u_2(t))>c_0 \left(1+t-(k-1)T\right)^{-q}\textrm{ and }l_0(2k)< k \right).
\end{array}
$$
Thus applying \eqref{Abstract_abis}, using $2(t-(k-1)T)>t$, it follows
\begin{equation}\label{Abs_b}
\Espace
\begin{array}{l}
\P\left( d_E(u_1(t),u_2(t))> 2^q c_0 \left(1+t\right)^{-q}\right)\leq
\P\left(l_0(2k)\geq k\right)+
2^q c_0 \left(1+t\right)^{-q}.
\end{array}
\end{equation}
In order to estimate $\P\left(l_0(2k)\geq k\right)$, we introduce the following notation
$$
l_0(\infty)=\lim \sup l_0.
$$
Taking into account \eqref{Abstract_ater}, we obtain that for $l<\infty$
$$\{l_0(\infty)=l\}=\{ l_0(k)=l, \textrm{ for any } k\geq l\}.$$
We deduce
\begin{equation}\label{Abs_c}
\P\left(l_0(2k)\geq k\right)\leq \P\left(l_0(\infty)\geq k\right).
\end{equation}
Taking into account \eqref{Abs_b},  \eqref{Abs_c}  and using a Chebyshev inequality,
it is sufficient to obtain that there exist $c_5>0$ such that
\begin{equation}\label{Abs_d}
\E \left(l_0(\infty)^{q}\right)\leq  c_5 \left(1+\Hcal(u_0^1)+\Hcal(u_0^2) \right).
\end{equation}

\subsection{Definition of a sequence of stopping times}

Let
$$
\tau=\min\left\{t\in T\N \,|\, \Hcal(u_1(t))+\Hcal(u_2(t))\leq R_0 \right\}.
$$
Then, the Lyapunov structure \eqref{Abstract_c} implies that there exist $\delta_0>0$
and $c_6>0$ such that
\begin{equation}\label{Abs_e}
\E \left(\exp\left(\delta_0 \tau \right)\right)\leq  c_6 \left(1+\Hcal(u_0^1)+\Hcal(u_0^2) \right).
\end{equation}
For a proof, see the proof of (1.56) at the end of the subsection 1.4 of \cite{ODASSO}.

We set
$$
\hat \sig =\min \left\{ k\in \N\backslash\{0\} \,|\, l_0(k)>1  \right\},\quad
\sig      =\hat \sig T.
$$
Clearly $\hat \sig = 1$ if the two solutions do not get coupled at time $0$ or $T$. Otherwise, they get coupled at $0$ or
$T$ and remain coupled until $\sig$.

From now, we fix $q_1\in (q,q_0-1)$. Let us assume for the moment that
there exists $p_\infty$ such that if $\Hcal_0\leq R_0$, then
\begin{equation}\label{Abs_f}
\Espace
\left\{
\begin{array}{l}
\E \left(\sig^{q_1} 1_{\sig<\infty}\right)\leq  K ,\\
\P\left( \sig=\infty \right)\geq p_\infty >0.
\end{array}
\right.
\end{equation}
The proof is given  at the end of this section.

Now we build a sequence of stopping times
$$
\Espace
\begin{array}{lcllcl}
\tau_0&=&\tau,\\
\hat \sig_{k+1} &=&\min \left\{ l\in \N\backslash\{0\} \,|\, lT>\tau_k \textrm { and }l_0(l) T>\tau_k+T  \right\},&
 \sig_{k+1}&=&\hat \sig_{k+1} \,T\\
\tau_{k+1}       &= & \sig_{k+1}+\tau o \theta_{\sig_{k+1}},
\end{array}
$$
 where $(\theta_t)_t$ is the shift operator.
The idea is the following. We wait the time $\tau_k$ to enter the ball of radius $R_0$. Then, if we do not start
 coupling at time $\tau_k$, we try to couple at time $\tau_k+T$. If we fail to start coupling at time $\tau_k$ or
$\tau_k+T$ we set $\sig_k=\tau_k+T$ else we set $\sig_k$ the time the coupling fails ($\sig_k=\infty$ if the
 coupling never fails). Then if $\sig_k<\infty$, we retry to couple after  entering
 in the ball of radius $R_0$. The fact that $R_0\geq d_0$
 implies that  $l_0(\tau_k)\in \{\tau_k,\infty\}$.

Note that we clearly have  $l_0(\tau_k)\in \{\tau_k,\infty\}$ and
$l_0(\sig_k)\in\{\sig_k,\infty\}$, and the  $l_0$--Markov property
implies the strong Markov property when conditioning
 with respect to $\F_{\tau_k}$ or $\F_{\sig_k}$.

We infer from the $l_0$--Markov property of $V$ that
$$
\sig_{k+1}=\tau_k+\sig o \theta_{\tau_k},
$$
which implies
$$
\tau_{k+1}=\tau_k+\rho o \theta_{\tau_k}, \; \textrm{ where } \rho=\sig + \tau o \theta_{\sig} .
$$

\subsection{Polynomial estimate on $\rho$}

We first establish that there exist $K_0$ such
that for any $V_0$ such that $\Hcal_0\leq R_0$
\begin{equation}\label{Abs_i}
\E_{V_0} \left(\rho^{q_1} 1_{\rho<\infty}\right) \leq  K_0.
\end{equation}

Notice that for any $V_0$ such that $\Hcal_0\leq R_0$,
\begin{equation}\label{Eq1.4.2}
\E_{V_0} \left(\rho^{q_1} 1_{\rho<\infty}\right) \leq c\left(
\E_{V_0} \left(\sig^{q_1 }1_{\sig<\infty}\right)+
\E \left(  ( \tau o \theta_{\sig} )^{q_1} 1_{\tau o \theta_{\sig}<\infty}1_{\sig<\infty}
\right)  \right).
\end{equation}
Applying the $l_0$--Markovian property and \eqref{Abs_e}, we obtain
$$
\E \left(  ( \tau o \theta_{\sig} )^{q_1} 1_{\tau o \theta_{\sig}<\infty}1_{\sig<\infty} | \F_{\sig}\right)
\leq  c_6 \left(1+\Hcal(u_1(\sig))+\Hcal(u_2(\sig)) \right)1_{\sig<\infty},
$$
which implies by applying the Lyapunov structure \eqref{Abstract_c}
\begin{equation}\label{Eq1.4.1}
\E \left(  ( \tau o \theta_{\sig} )^{q_1} 1_{\tau o \theta_{\sig}<\infty}\right)\leq  c_6 (1+2K' (R_0+\E(\sigma1_{\sig<\infty})).
\end{equation}
Applying \eqref{Abs_f} and \eqref{Eq1.4.1} to \eqref{Eq1.4.2}, we obtain \eqref{Abs_i}.

\subsection{Conclusion}

Applying a convexity inequality, we obtain
$$
\E \left(\tau_k^{q_1}1_{\tau_k<\infty}\right)\leq  (k+1)^{(q_1-1)^+}\left(\E\tau^{q_1}+
\sum_{n=0}^{k-1}\E (\rho o \theta_{\tau_n})^{q_1}1_{\rho o \theta_{\tau_n}<\infty}\right).
$$
Combining the $l_0$--Markov property, \eqref{Abs_e} and
\eqref{Abs_i} gives
\begin{equation}\label{abs_c}
\E \left(\tau_k^{q_1}1_{\tau_k<\infty}\right)\leq  C (k+1)^{1\vee q_1}\left(1+\Hcal(u_0^1)+\Hcal(u_0^2) \right).
\end{equation}

Now, we are able to estimate $\E \left(l_0(\infty)^{q}\right)$
$$
\E \left(l_0(\infty)^{q}\right)\leq c\left(1+\sum_{n=0}^\infty \E \left(\tau_n^{q}1_{\tau_n<\infty}1_{k_0=n}\right) \right),
$$
where
$$k_0=\inf\{k\in \N\,|\, \sig_{k+1}=\infty\}.$$
Then, applying an Holder inequality, we obtain
$$
\E \left(l_0(\infty)\right)^{q}\leq
c\left(1+ \sum_{n=0}^\infty \E \left(\tau_n^{pq}1_{\tau_n<\infty}\right)^{\frac{1}{p}}
\left(\P\left( k_0=n \right)\right)^\frac{1}{p'} \right).
$$
Using the second inequality of \eqref{Abs_f} and  $\tau<\infty$, we obtain from the $l_0$--Markov property that
\begin{equation}\label{Abs_k}
\P\left( k_0>n \right)\leq \left(1-p_\infty\right)^n.
\end{equation}
It follows that  $k_0<\infty$ almost surely and that
$$
l_0(\infty)\in\{\tau_{k_0},\tau_{k_0}+1\}.
$$
Therefore $l_0(\infty)<\infty$ almost surely and applying \eqref{abs_c}, we obtain that if $pq=q_1$
$$
\E \left(l_0(\infty)\right)^{q}\leq
C\left(\sum_{n=0}^\infty (n+1)^{\frac{1}{p}\vee q}
(1-p_\infty)^\frac{n}{p'}\right)\left(1+\Hcal(u_0^1)+\Hcal(u_0^2) \right).
$$
Thus \eqref{Abs_d} is established and we can conclude.

\subsection{Proof of \eqref{Abs_f}}

Now we establish \eqref{Abs_f}. There are two cases. The first case is $l_0(0)=0$. Then, applying \eqref{Abstract_a},
we obtain that
$$
\P\left( \sig=\infty \right)\geq\Pi_{k=0}^\infty\P\left( l_0(k+1)=0|l_0(k)=0 \right)\geq \Pi_{k=0}^\infty p_k.
$$
The second case is $l_0(0)=\infty$. Then
$$
\P\left( \sig=\infty \right)\geq\P\left( l_0(1)=1 \right)\Pi_{k=1}^\infty\P\left( l_0(k+1)=1|l_0(k)=1 \right).
$$
Since $\Hcal_0\leq R_0$, then applying \eqref{Abstract_a} and \eqref{Abstract_b}
$$
\P\left( \sig=\infty \right)\geq \Pi_{k=-1}^\infty p_k.
$$
Since $p_k>0$ and $1-p_k$ decreases to $0$ faster than $k^{-q_0}$ with $q_0>1$, then the product converges and
 in the two cases
\begin{equation}\label{Abs_l}
\P\left( \sig=\infty \right)\geq p_\infty =\Pi_{k=-1}^\infty p_k>0.
\end{equation}
Notice that \eqref{Abstract_a} implies
$$
\begin{array}{rcl}
\P\left( \sig = n \right) & \leq &
 \P\left( l_0(n+1)\not= n\,|\,l_0(n)=0\right)+ \P\left( l_0(n+1)\not= n\,|\,l_0(n)=1\right),\\
&\leq & 2c_1 \left(1+(n-1)T \right)^{-q_0},
\end{array}
$$
which gives the first inequality of \eqref{Abs_f} and allows to conclude because $q_1<q_0-1$.


\section{Proof of the  a priori estimates}

As already mentioned, the various computations made in this section are not rigorous. Indeed, the
solutions are not smooth enough to apply Ito formula. A suitable approximation could be used to justify
the result rigorously.

{\bf Ito Formula for $\abs{u}^6$}

Applying Ito Formula to $\abs{u}^6$, we obtain
$$
d \abs{u}^6 + 6\alpha\abs{u}^6 dt =6 \abs{u}^4\left(u,bdW\right)+12 \abs{u}^{2} \abs{b^* u}^2 dt+3B_0\abs{u}^4dt.
$$
Since $b^*$ is a bounded operator on $L^2(0,1)$,
$$
12 \abs{u}^{2} \abs{b^* u}^2\leq 12 B_0 \abs{u}^4.
$$
We deduce
$$
12 \abs{u}^{2} \abs{b^* u}^2 + 3 B_0 \abs{u}^4\leq \alpha \abs{u}^6 + C,
$$
and
\begin{equation}\label{1.3.1}
d \abs{u}^6 + 5\alpha\abs{u}^6 dt \leq 6 \abs{u}^4\left(u,bdW\right)+Cdt.
\end{equation}

{\bf Ito Formula for $\Hcal$}

Applying Ito Formula to $\Hcal_*$, we obtain
\begin{equation}\label{1.3.2}
d\Hcal_*(u)+\alpha\left( \norm{u}^2-\abs{u}_4^4 \right)dt = dM_*+B_1dt +I_*dt,
\end{equation}
where
$$
\Espace
\begin{array}{lcl}
dM_*&=& \left( A u -\abs{u}^2u, bdW  \right), \\
I_* &=& -\sum_{n=1}^\infty b_n^2 \int_{[0,1]}
\left(2\Re(u(t,x)\overline{ e_n}(x))^2+\abs{e_n(x)}^2\abs{u(t,x)}^2  \right)dx.
\end{array}
$$
Note that, since $b^*A$ is a bounded operator from $L^2(0,1)$ to $H^1_0(0,1)$, $M_*$ is well defined.
Recalling that $\abs{e_n}_\infty=1$, we obtain
$$
I_*\leq 3 B_0 \abs{u}^2\leq \alpha c_0 \abs{u}^6 + C.
$$
Recalling that $
\abs{\cdot}^4_4\leq \frac{1}{4}\norm{\cdot}^2+c_0\abs{\cdot}^6
$,
 we infer from \eqref{1.3.1}, \eqref{1.3.2} and the last inequality that
\begin{equation}\label{1.3.3}
d\Hcal(u)+\frac{3}{2}\alpha \Hcal(u) dt \leq dM_1+C_1dt,
\end{equation}
where
$$
dM_1=dM_*+6 c_0 \abs{u}^4\left(u,bdW\right).
$$

{\bf Ito Formula for $\Hcal^k$}

Applying Ito Formula to $\Hcal^k$ for $k\in \N\backslash\{0\}$, we obtain similarly as above
\begin{equation}\label{1.3.4}
d\Hcal(u)^k+\frac{3}{2}\alpha k \Hcal(u)^k dt \leq dM_k+k\Hcal(u)^{k-1}C_1dt+\frac{k(k-1)}{2}\Hcal(u)^{k-2}d\left<M_1\right>,
\end{equation}
where
$$
dM_k=k\Hcal(u)^{k-1}dM_1.
$$
Note that, since $b^*A$ is a bounded operator from $L^2(0,1)$ to $H^1_0(0,1)$ and $b^*$ is bounded
from $L^1(0,1)$ to $L^2(0,1)$,
$$
\abs{b^*\left(Au-\abs{u}^2 u\right)}^2\leq 4B_1\|u\|^2+cB_1\abs{u}_3^6,
$$
it follows from a Gagliardo-Nirenberg inequality
$$
\abs{b^*\left(Au-\abs{u}^2 u\right)}^2\leq cB_1\left( \norm{u}^2+\abs{u}^{10} \right).
$$
Now, we write
$$
d\left<M_1\right>\leq 2\abs{b^*\left(Au-\abs{u}^2 u\right)}^2+72 c_0^2 \abs{u}^8\abs{b^*u}^2,
$$
and deduce that
$$
d\left<M_1\right>\leq cB_1\left( \norm{u}^2+\abs{u}^{10} \right),
$$
and
\begin{equation}\label{1.3.5}
d\left<M_1\right>\leq cB_1\left(1+\Hcal(u)^{\frac{5}{3}}\right).
\end{equation}

Gathering \eqref{1.3.4} and \eqref{1.3.5} and using once more an arithmetico-geometric
inequality, we obtain
\begin{equation}\label{1.3.6}
d\Hcal(u)^k+\alpha k \Hcal(u)^k dt \leq dM_k+C_k''dt.
\end{equation}

{\bf Proof of Lemma \ref{lem_Lyapounov_norm}}

Multiplying  \eqref{1.3.6} by $e^{\alpha kt}$ yields
$$
d\left(e^{\alpha kt}\Hcal(u)^k\right)\leq e^{\alpha kt}dM_k+e^{\alpha kt}C_k''dt.
$$
By integration we obtain
$$
e^{\alpha kt} \mathcal H(u(t))^k\leq  \Hcal(u_0)^k + \int_0^t e^{\alpha ks}  dM_k(s)+
\frac{C_k''}{\alpha k} e^{\alpha kt}
$$
and
$$
\mathcal H(u(t))^k\leq  \Hcal(u_0)^k e^{-\alpha k t}+ \int_0^t e^{-\alpha k(t-s)}  dM_k(s)+
\frac{C_k''}{\alpha k},
$$
which yields, by taking the expectation, the first inequality of Lemma \ref{lem_Lyapounov_norm}.

\noindent Let $M>0$ and $\tau\leq M$ be a bounded stopping time. Then, integrating \eqref{1.3.6} between $0$ and $\tau$and taking the expectation yields
$$
\E\left(\Hcal(u(\tau))^k\right)\le \Hcal(u_0)^k +C"_k\E(\tau),
$$
which is  the second inequality of Lemma \ref{lem_Lyapounov_norm} for bounded stopping times.

\noindent Assume now that $\tau$ is a general stopping time. We consider the second inequality of Lemma
\ref{lem_Lyapounov_norm}
for the stopping time $\tau\wedge M$ and we take the limit when $M\to \infty$.
The second inequality of Lemma \ref{lem_Lyapounov_norm} for $\tau$ follows from Fatou Lemma.

{\bf Proof of Proposition \ref{Prop_majoration_energie_sous_lineaire_Hun}}

We first note that
$$
d\left< M_k \right>=k^2\Hcal(u)^{2(k-1)}  d\left< M_1 \right>,
$$
so that,  taking into account \eqref{1.3.5},
\begin{equation}\label{1.3.7}
d\left< M_k \right>\leq c_k\left(1+ \Hcal(u)^{2k}\right)ds
\end{equation}

Taking the expectation of \eqref{1.3.6}, we obtain for any $k\ge 1$
\begin{equation}\label{1.3.9}
\E \int_0^t \Hcal(u(s))^k dt \leq C_k\left( \Hcal(u_0)^k + t\right).
\end{equation}
Hence, for any $p\ge 1$,
\begin{equation}\label{1.3.10}
\E \left< M_k \right>^p(t) \leq C_{k,p}\left( \Hcal(u_0)^{2kp} + t^p\right).
\end{equation}

Applying the maximal martingale inequality and taking into account \eqref{1.3.6}, we infer from \eqref{1.3.10} the first
 inequality of Proposition
\ref{Prop_majoration_energie_sous_lineaire_Hun}.

Applying the maximal martingale inequality on $[n,n+1]$, $n\geq 0$, we have
$$
\P\left( \sup_{[n,n+1]} M_k > a+\Hcal(u_0)^{2k}+n+1  \right)
\leq c_p \frac{\E \left< M_k \right>^{p+1}(n+1)}{ \left(a+\Hcal(u_0)^{2k}+n+1 \right)^{2p+2}}.
$$
It follows from \eqref{1.3.10} that
\begin{equation}\label{1.3.14}
\P\left( \sup_{[n,n+1]} M_k > a+\Hcal(u_0)^{2k}+n+1  \right)
\leq \frac{c_p C_{k,p+1}'}  {\left(a+\Hcal(u_0)^{2k}+n+1 \right)^{p+1}}.
\end{equation}
Now, summing \eqref{1.3.14} over $n\geq T$, for
$T$ integer, we obtain that for any $(p,k)\in (\N\backslash\{0\})^2$ there exists $K_{k,p}$ such that
\begin{equation}
 \label{1.3.12}
 \P\left( \sup_{t\in [T,\infty)} M_k(t) >  1+a+\Hcal(u_0)^{2k}+t \right)
\leq  K_{k,p}\left( a+T \right)^{-p} ,\; T>0.
\end{equation}
Taking into account \eqref{1.3.6}, this implies the second inequality of Proposition
\ref{Prop_majoration_energie_sous_lineaire_Hun}.

{\bf Proof of Proposition \ref{Prop_petit}}

Combining Lemma \ref{lem_Lyapounov_norm} applied to $\tau=t$ and  Chebyshev's inequality, we obtain
\begin{Lemma}\label{lem_Lyapounov_proba}
Let $(u_i,W_i)_{i=1,2}$ be a couple of solutions of \eqref{Eqbase}, \eqref{Eqinitial} such that $W_1$ and $W_2$ are two
 cylindrical Wiener process on $L^2([0,1])$. If $R_0 \geq  \left(\sum_{i=1}^2\Hcal(u_0^i)\right) \vee C_1 $, then
$$
\P\left ( \mathcal H(u_1(t))+\Hcal(u_2(t)) \geq  4 C_1  \right ) \leq \frac{1}{2},
$$
provided
$t \geq  \theta_1(R_0)= \frac{1}{\alpha} \ln \frac{R_0}{C_1}$.
\end{Lemma}

It follows from Lemma \ref{lem_Lyapounov_proba} that it is sufficient to establish Proposition \ref{Prop_petit}
 for $R_0=4C_1$ and $t=T_{-1}(R_0,R_1)$ (instead of $t\geq T_{-1}(R_0,R_1)$). From now on, we only consider the case
 $R_0=4C_1$.

Let $T,\delta>0$. Applying Chebyshev inequality, we obtain ${N_{-2}}=N_{-2}(T,\delta)\in \N$ such that
$$
\P\left( \sup_{t\in [0,T]}\norm{bQ_{N_{-2}}W(t)}_3 > \frac{\delta}{2} \right)
\leq \frac{2}{\delta}\sum_{n> N_2}\mu_n^3 b_n^2 \leq \frac{1}{2}.
$$
Moreover $P_{N_{-2}}W$ is a finite dimensional brownian motion and
it is classical that
$$
\pi_{-3}(T,\delta,N_{-2})=
\P\left( \sup_{t\in [0,T]}\abs{P_{N_{-2}}W(t)} \leq \frac{\delta}{2}\norm{b}_{\mathcal L_2(L_2(D),H^3(D))}^{-1} \right) >0.
$$
Writing
$$
\P\left( \sup_{t\in [0,T]}\norm{bW(t)}_3 \leq \delta \right)\geq
 \P\left( \sup_{t\in [0,T]}\norm{bQ_{N_{-2}}W(t)}_3 \leq \frac{\delta}{2} \right)
\pi_{-3}(T,\delta,N_{-2}),
$$
it follows
\begin{equation}\label{petit1}
\pi_{-2}(T,\delta)=\P\left( \sup_{t\in [0,T]}\norm{bW(t)}_3 \leq \delta \right)>0.
\end{equation}

It thus suffices to prove that there exists $T_{-1}(R_1),\delta_{-1}(R_1)>0$ such that
\begin{equation}\label{petit2}
\left\{ \sup_{t\in [0,T_{-1}]}\norm{bW(t)}_3 \leq \delta_{-1} \right\}
\subset
\left\{ \Hcal(u(T_{-1},u_0))\leq \frac{1}{2}R_1 \right\},
\end{equation}
provided $\Hcal(u_0)\leq R_0$.

{\bf Proof of \eqref{petit2}}

Let us set
$$
v=u(\cdot,u_0)-bW,
$$
then
\begin{equation}\label{petit3}
\frac{dv}{dt}+\alpha v +\i A v -\i \Bd{bW+v}=-\left(\alpha +\i A\right)bW.
\end{equation}
Taking the scalar product between \eqref{petit3} and $2v$, we obtain
$$
\frac{d\abs{v}^2}{dt}+2\alpha \abs{v}^2 = 2\left(v,\i \Bd{bW+v}-\left(\alpha +\i A\right)bW\right).
$$
Since
$$
\left(v,\i \abs{bW+v}^2v\right)=0,
$$
applying H\"older inequalities and Sobolev Embedding $H^1(D)\subset L^\infty(0,1)$, we deduce
$$
\frac{d\abs{v}^2}{dt}+2\alpha \abs{v}^2 \leq c\norm{bW}_3\left(1+\norm{bW}_3^2\right)\left(1+\norm{v}^3\right).
$$
Applying Ito Formula to $\abs{v}^6$, we deduce
\begin{equation}\label{petit4}
\frac{d\abs{v}^6}{dt}+6\alpha \abs{v}^6 \leq c\norm{bW}_3\left(1+\norm{bW}_3^2\right)\left(1+\norm{v}^9\right).
\end{equation}
Taking the scalar product between \eqref{petit3} and $Av-\B v$, we obtain
$$
\frac{d\Hcal_*(v)}{dt}+\alpha\norm{v}^2=-\left(Av-\B v,\left(\alpha +\i A\right)bW\right)
+\alpha\left( \Bd{v+bW},v\right).
$$
Since
$$
I_1=\alpha\left(\left( \Bd{v+bW},v\right)-\abs{v}_4^4\right)=\alpha\left( \Bd{v+bW}-\B v,v\right),
$$ we obtain
\begin{equation}\label{petit5}
\frac{d\Hcal_*(v)}{dt}+\alpha\left(\norm{v}^2-\abs{v}_4^4\right)=I_1+I_2,
\end{equation}
where
$$
I_2=-\left(Av-\B v,\left(\alpha +\i A\right)bW\right)
.
$$
Recalling that for any $z,h\in \C^2$
$$
\abs{\Bd{z+h}-\B z}\leq \abs{h}\left(\abs{z}^2+\abs{h}^2\right),
$$
and applying H\"older inequality and the Sobolev Embedding  $H^1(D)\subset L^\infty(0,1)$ , we obtain
$$
I=I_1+I_2\leq c\norm{bW}_3\left(1+\norm{v}^3\right)\left(1+\norm{bW}_3^2\right).
$$
It follows from \eqref{petit4}, \eqref{petit5} and the last inequality  that
\begin{equation}\label{petit6}
\frac{d\Hcal(v)}{dt}+2\alpha\Hcal(v)\leq c\norm{bW}_3\left(1+\norm{bW}_3^2\right)\left(1+ \Hcal(v)^5 \right).
\end{equation}

Let $T,\delta,M>0$ and assume that
$$
 \sup_{t\in [0,T]}\norm{bW(t)}_3 \leq \delta.
$$
We set
$$
\tau=\inf\left\{t\in [0,T]\;|\; \Hcal(v)\leq 3R_0  \right\}.
$$
Integrating \eqref{petit6}, we obtain
\begin{equation}\label{petit7}
\Hcal(v(t))\leq e^{-2\alpha t} R_0+\frac{c}{2\alpha} \delta\left(1+\delta^2\right)\left(1+ R_0^5 \right),
\end{equation}
provided $t\leq \tau$.

Now we choose $\delta\leq \delta_{-2}(R_1')>0$ such that
$$
\frac{c}{2\alpha} \delta\left(1+\delta^2\right)\left(1+ R_0^5 \right)\leq R_1'\wedge R_0.
$$
It follows from \eqref{petit7} that
$$
\tau=T,
$$
and that
$$
\Hcal(v(T))\leq 2 R_1',
$$
provided
$$
\quad T\geq \frac{1}{2\alpha}\ln \left( \frac{R_1'}{R_0}\right).
$$
In order to conclude, we remark that
$$
\Hcal(u(T))\leq c\left(\Hcal(bW(T))+ \Hcal(v(T))\right)\leq c(\delta^2(1+\delta^4)+R_1').
$$
Then, choosing $\delta$ and $R_1'$ sufficiently small, we obtain \eqref{petit2}.

\section{Proof of the Foias-Prodi estimates}

The aim of this section is to establish Proposition \ref{Prop_Foias_Prodi_Hun}.

{\bf $L^2$ estimates}

\noindent Taking into account \eqref{1.2.1}, we deduce that the difference of the two solutions $r=u_1-u_2$ satisfies the equation
\begin{equation}\label{1.4.1}
\LS{r}=\i Q_N\left( \abs{u_1}^2u_1-\abs{u_2}^2u_2  \right).
\end{equation}
Applying Ito Formula to $\abs{r}^2$, we obtain
$$
\LSabs{r}=2\left( \i r, \abs{u_2}^2u_2- \abs{u_1}^2u_1 \right).
$$
Since
$$
\abs{\abs{u_2}^2u_2- \abs{u_1}^2u_1}\leq c\left( \sum_{i=1}^2 \abs{u_i}^2 \right)\abs{r},
$$
it follows
$$
\LSabs{r}\leq c\int_{[0,1]} \left(\sum_{i=1}^2 \abs{u_i}^2\right) \abs{r}^2 dx.
$$
Using the Sobolev Embedding $H^1(0,1)\subset L^\infty(0,1)$, we obtain
\begin{equation}\label{1.3.3tierce}
\LSabs{r}\leq c \left(\sum_{i=1}^2 \Hcal(u_i)\right) \abs{r}^2.
\end{equation}
We deduce as in the proof of \eqref{1.3.3}
$$
d\Hcal(u_i)+\frac{3}{2}\alpha \Hcal(u_i) dt \leq dM_1^i+C_1dt+1_{i=1}(G,h)dt,
$$
where
$$
\Espace\left\{
\begin{array}{rcl}
dM_1^i&=&\left( A u_i -\abs{u_i}^2u_i, bdW_i  \right) +6 c_0 \abs{u_i}^4\left(u_i,bdW_i\right),\\
G&=& A u_1 -\abs{u_1}^2u_1+6 c_0 \abs{u_1}_{L^2}^4u_1.
\end{array}\right.
$$
It follows from Sobolev Embeddings and H\"older inequalities that
$$
\norm{G}_{-1}\leq c\left(1+\Hcal(u_1)\right)^\frac{5}{6}.
$$
Hence we deduce from \eqref{1.2.2} that
$$
(G,h)\leq c\left(1+\Hcal(u_1)+\Hcal(u_2)\right)^4.
$$
Taking into account \eqref{1.3.3tierce}, it follows
\begin{equation}\label{1.4.2}
  dZ_1+2\alpha Z_1 dt \leq c \left(1+\sum_{i=1}^2 \Hcal(u_i)^4\right) \abs{r}^2dt+\abs{r}^2 dM_\#,
\end{equation}
where
$$
Z_1=\left(\sum_{i=1}^2\Hcal(u_i)\right)\abs{r}^2
$$
and
$$
dM_\#=dM_1^1+dM_1^2.
$$

{\bf Ito Formula for $J$}

Now we rewrite \eqref{1.4.1} in the form
\begin{equation}\label{1.4.6}
\LS{r}=-\i\frac{1}{2}Q_N\left( (\abs{u_1}^2+\abs{u_2}^2)r+\Re\left((u_1+u_2)\bar r  \right) (u_1+u_2) \right).
\end{equation}

Applying Ito Formula to $J_*(u_1,u_2,r)$, we obtain
\begin{equation}\label{1.4.7}
\Espace
\begin{array}{r}
dJ_*+2\alpha J_*dt=g(u_1,u_2,r)dt+g(u_2,u_1,r)dt +\psi(u_1,u_2,r)(bdW_1)\\
 +\psi(u_1,u_2,r)(h(t))dt+\psi(u_2,u_1,r)(bdW_2)+I_1(r)dt+dI_2(r,dt),
\end{array}
\end{equation}
where
$$
\Espace
\begin{array}{lcl}
g(u_1,u_2,r) & = & \left\{
\begin{array}{ll}
 &2 \int_{[0,1]} \left(\Re\left(\bar u_1(\alpha u_1+\i A u_1 -\i \B{u_1} )\right)\abs{r}^2  \right) dx \\
+ &2 \int_{[0,1]} \Re\left(\bar r( u_1+u_2 )\right)\Re\left( \bar r(\alpha u_1-\i A u_1 +\i \B{u_1} ) \right) dx
\end{array}
\right.,\\
\psi(u_1,u_2,r)(h) & = &
 2 \int_{[0,1]} \left(\Re\left(\bar u_1 h\right)\abs{r}^2  \right) dx
+ 2 \int_{[0,1]} \Re\left(\bar r( u_1+u_2 )\right)\Re\left( \bar r h \right) dx
,\\
I_1(r)&=& -\sum_{n=1}^\infty b_n^2\int_{[0,1]} \left( \abs{e_n}^2\abs{r}^2+\Re(e_n\bar r)^2 \right)dx,\\
dI_2(r,t)&=& -\sum_{p,q=1}^\infty b_p b_q
\left(\left(\int_{[0,1]} \Re\left( e_p \bar r \right)\Re\left(e_q\bar r\right)dx\right)\; d\left< (W_1,e_p),(W_2,e_q) \right>\right).
\end{array}
$$
Applying an integration by part to $Au_1$, H\"older inequality and the Sobolev Embedding
$H^{\frac{3}{4}}(0,1)\subset L^\infty(0,1)$, we obtain
\begin{equation}\label{1.4.8}
g(u_1,u_2,r)\leq \left(1+\sum_{i=1}^2\norm{u_i}^6  \right)\norm{r}\norm{r}_\frac{3}{4}.
\end{equation}
We deduce from  H\"older inequality that
$$
\psi(u_1,u_2,r)(h(t))\leq\left(\sum_{i=1}^2 \abs{u_i}_\infty\right)\abs{h(t)}\abs{r}^2_4.
$$
Taking into account \eqref{1.2.2} and applying the Sobolev Embeddings $H^1(0,1)\subset L^\infty(0,1)$ and
$H^\frac{1}{2}(0,1)\subset L^4(0,1)$, we obtain
\begin{equation}\label{1.4.9}
\psi(u_1,u_2,r)(h(t))\leq c\kappa_0\left(1+\sum_{i=1}^2 \Hcal(u_i)^\frac{7}{2}\right)\norm{r}_\frac{1}{2}^2.
\end{equation}
Recalling that $\abs{e_n}_\infty=1$, we obtain
\begin{equation}\label{1.4.10}
I_1(r)\leq 3B_0\abs{r}^2.
\end{equation}
Note that we have no information on the law of the couple $(W_1,W_2)$.
Hence, we cannot compute $d\left< (W_1,e_p),(W_2,e_q) \right>$. However we know that
$$
d\abs{\left< (W_1,e_p),(W_2,e_q) \right>}\leq dt.
$$
Hence
$$
d\abs{I_2(r,t)}=
\left(\int_{[0,1]} \Re\left( \sum_{n=1}^\infty \left(b_n e_n\right) \bar r \right)^2dx\right)\; dt.
$$
Applying  the following Schwartz  inequality
$$
\left(\sum_{n=1}^\infty b_n \right)^2\leq  \left(\sum_{n=1}^\infty \mu_n b_n^2\right)
 \left(\sum_{n=1}^\infty \frac{1}{\mu_n}\right)\leq cB_1,
$$
we deduce from $\abs{e_n}_\infty=1$ that
\begin{equation}\label{1.4.11}
d\abs{I_2(r,t)}\leq cB_1\abs{r}^2dt.
\end{equation}

Combining \eqref{1.4.7}, \eqref{1.4.8}, \eqref{1.4.9}, \eqref{1.4.10}, and \eqref{1.4.11}  , we obtain
\begin{equation}\label{1.4.12}
dJ_*+2\alpha J_*dt\leq c \left(1+\sum_{i=1}^2 \Hcal(u_i)^{4}\right)\norm{r}\norm{r}_{\frac{3}{4}}dt+dM_{\#\#},
\end{equation}
where
$$
dM_{\#\#}=\left(\psi(u_1,u_2,r)(bdW_1)+\psi(u_2,u_1,r)(bdW_2)\right).
$$
Summing \eqref{1.4.2} and \eqref{1.4.12}, we obtain
\begin{equation}\label{1.4.13}
dJ+2\alpha Jdt\leq c \left(1+\sum_{i=1}^2 \Hcal(u_i)^{4}\right)\norm{r}\norm{r}_{\frac{3}{4}}dt+dM,
\end{equation}
where
$$
dM=dM_{\#\#}+c_1\abs{r}^2dM_\#.
$$

{\bf Conclusion}

Since $\norm{r}_{\frac{3}{4}}\leq \mu_{N+1}^{-\frac{1}{8}}\norm{r}$ then there exists $\Lambda>0$ such that
\begin{equation}\label{1.4.3}
dJ+\left(2\alpha  - \frac{\Lambda}{\mu_{N+1}^{\frac{1}{8}}} l(u_1,u_2)\right) J dt\leq dM.
\end{equation}

Multiplying \eqref{1.4.3} by $\exp\left(2\alpha s-\Lambda \mu_{N+1}^{-\frac{1}{8}} \int_0^s l(u_1(s'),u_2(s'))ds' \right)$,
 we obtain that
\begin{equation}\label{1.4.5}
J_{FP}^N(t\wedge\tau)\leq
\int_0^{t\wedge\tau}
\exp\left(\frac{3}{2}\alpha s-\Lambda \mu_{N+1}^{-\frac{1}{8}} \int_0^s l(u_1(s'),u_2(s'))ds' \right)
dM(s).
\end{equation}
Fatou Lemma allows to conclude.

\vspace{.5cm}

{\bf Acknowledgments: } We are very grateful to the referee for his many suggestions, which have led
to many improvements in this article . We also would like to thank A. Shirikyan for many fruitful discussions.






\end{document}